\documentclass[11pt,reqno]{amsart} 

\usepackage{amsthm, amsfonts, amssymb, amsmath, enumerate, bbm,paralist,enumitem,xr,pkgfile,float,caption,subcaption,qtree,newfloat}
\RequirePackage[numbers]{natbib}
\RequirePackage{graphicx}

\newtheorem{theorem}{Theorem}
\newtheorem{lemma}{Lemma}

\newtheorem{corollary}{Corollary}

\newcommand{\bmp}{{\hat{\beta}_{\textrm{MPL}}}} 
\newcommand{\bml}{{\hat{\beta}_{\textrm{ML}}}}

\newcommand{\er}{{Erd\H{o}s-R\'enyi}}
\newcommand{\sa}{{\{-1,1\}^n}}
\newcommand{\p}{{\mathbb{P}}}

\newcommand{\mn}{{\mathcal{M}_n}}

\newcommand{\e}{{\mathbb{E}}}
\newcommand{\bs}{{\bm X}}
\newcommand{\os}{{\overline{X}_n}}
\newcommand{\oss}{{\overline{x}_n}}
\newcommand{\ml}{{\mathscr{L}_n(\bm x)}}
\newcommand{\bt}{{\bm x}}

\newcommand{\ern}{{Erd\H{o}s-R\'enyi}}

\newcommand{\um}{{\underline{m}}}
\newcommand{\hbp}{{H_{\beta,p}}}
\newcommand{\hbop}{{H_{\beta_0,p}}}
\allowdisplaybreaks

\newcounter{myalgctr}

\numberwithin{myalgctr}{section}

\date{} 
\title[Efficient Estimation in Tensor Ising Models]{Efficient Estimation in Tensor Ising Models}

\author[Mukherjee]{Somabha Mukherjee} 
\address{Department of Statistics and Data Science, National University of Singapore {\tt somabha@nus.edu.sg}}

\author[Son]{Jaesung Son} 
\address{Department of Statistics, Columbia University {\tt js4638@columbia.edu}}

\author[Ghosh]{Swarnadip Ghosh} 
\address{Department of Statistics, Stanford University {\tt raswa281@stanford.edu}}

\author[Mukherjee]{Sourav Mukherjee} 
\address{Department of Statistics, University of Florida {\tt souravmukherjee@ufl.edu}}
\begin{document}
	\maketitle
	\begin{abstract}
		The tensor Ising model is a discrete exponential family used for modeling binary data on networks with not just pairwise, but higher-order dependencies. A particularly important class of tensor Ising models are the tensor Curie-Weiss models, where all tuples of nodes of a particular order interact with the same intensity. The maximum likelihood estimator (MLE) is not explicit in this model, due to the presence of an intractable normalizing constant in the likelihood, and a computationally efficient alternative is to use the maximum pseudolikelihood estimator (MPLE). In this paper, we show that the MPLE is in fact as efficient as the MLE (in the Bahadur sense) in the $2$-spin model, and for all values of the null parameter above $\log 2$ in higher-order tensor models. Even if the null parameter happens to lie within the very small window between the threshold and $\log 2$, they are equally efficient unless the alternative parameter is large. Therefore, not only is the MPLE computationally preferable to the MLE, but also theoretically as efficient as the MLE over most of the parameter space. Our results extend to the more general class of Erd\H{o}s-R\'enyi hypergraph Ising models, under slight sparsities too.
		
	\end{abstract}
	
	\keywords
	
	
	
	
	
	\section{Introduction}
	
	With the ever increasing demand for modeling dependent network data in modern statistics, there has been a noticeable rise in the necessity for introducing appropriate statistical frameworks for modeling dependent data in the recent past. One such useful and mathematically tractable model which was originally coined by physicists for describing magnetic spins of particles, and later used by statisticians for modeling dependent binary data, is the Ising model \cite{ising}. It has found immense applications in diverse places such as image processing \cite{geman_graffinge}, neural networks \cite{neural}, spatial statistics \cite{spatial}, disease mapping in epidemiology \cite{disease}, structure detection \cite{neykov1} and property testing \cite{neykov2}.
	
	The Ising model is a discrete exponential family on the set of all binary tuples of a fixed length, with sufficient statistic given by a quadratic form, designed to capture pairwise dependence between the binary variables, arising from an underlying network structure. However, in most real-life scenarios, pairwise interactions are not enough to capture all the complex dependencies in a network data. For example, the behavior of an individual in a peer group depends not just on pairwise interactions, but is a more complex function of higher order interactions with colleagues. Similarly, in physics, it is known that the atoms on a crystal surface do not just interact in pairs, but in triangles, quadruples and higher order tuples. A useful framework for capturing such higher order dependencies is the $p$-tensor Ising model \cite{jaesung1}, where the quadratic interaction term in the sufficient statistic is replaced by a multilinear polynomial of degree $p\ge 2$. Although constructing consistent estimates of the natural parameter in general $p$-tensor Ising models is possible \cite{jaesung1}, more exact inferential tasks such as constructing confidence intervals and hypothesis testing is not possible, unless one imposes additional constraints on the underlying network structure. One such useful structural assumption is that all $p$-tuples of nodes in the underlying network interact, and that too with the same intensity. The corresponding model is called the $p$-tensor Curie-Weiss model \cite{jaesungcw}, which is a discrete exponential family on the hypercube $\{-1,1\}^n$, with probability mass function given by:
	\begin{equation}\label{cwmodel}
		\p_{\beta,p}(\bm x) := \frac{\exp\left\{\beta n^{1-p} \sum_{1\le i_1,\ldots,i_p \le n} x_{i_1}\ldots x_{i_p} \right\}}{2^n Z_n(\beta,p)}\quad\textrm{for} ~\bm x \in \{-1,1\}^n~.
	\end{equation}
	Here $Z_n(\beta,p)$ is a normalizing constant required to ensure that $\sum_{\bm x \in \{-1,1\}^n} \p_{\beta,p}(\bm x) = 1$, and $\beta \ge 0$. It is precisely this inexplicit normalizing constant $Z_n(\beta,p)$, that hinders estimation of the parameter $\beta$ using the maximum likelihood (ML) approach. Although the ML estimator can still be computed in $O(n)$ time in the Curie-Weiss model (since the probability mass function \eqref{cwmodel} is actually a function of the sum $\sum_{i=1}^n x_i$ which can take $2n+1$ values), in even slightly more general models, for example the Erd\H{o}s-R\'enyi Ising model \eqref{ermodel}, this is not true, and the ML estimator is computationally infeasible.
	
	An extremely useful approach in the literature to circumvent this issue, is the concept of \textit{maximum pseudolikelihood} (MPL) estimation, which was introduced by Besag in the context of spatial stochastic data with both lattice and non-lattice interactions \cite{besag_lattice,besag_nl}, and is based on computing explicit conditional distributions.  This approach was later applied by Chatterjee to parameter estimation in Ising models (for $p=2$) under general network interactions \cite{chspin}. To elaborate, the MPL estimate is obtained by maximizing the \textit{pseudolikelihood function}:
	
	$$\hat{\beta}_{\textrm{MPL}} := \arg \max_{\beta \in \mathbb{R}} \prod_{i=1}^n \p_{\beta,p} \left(X_i | (X_j)_{j\ne i}\right)~,$$
	where $\bm X = (X_1,\ldots,X_n)$ is simulated from the model \eqref{cwmodel}. It should be quite clear that the conditional distributions in the expression of the pseudolikelihood function makes it free of the inexplicit normalizing constant $Z_n(\beta,p)$, thereby making the MPL estimator explicit and computationally feasible. In fact, methods as simple as a grid search can be applied to compute the MPL estimator.
	
	A natural question is thus, to what extent does this computationally feasible MPL approach inherit the desirable theoretical properties of the ML approach? Quite surprisingly, in spite of being just a proxy for the exact ML estimator, the MPL estimator in fact satisfies almost all the theoretical guarantees of the former. We list some of these properties below:
	\begin{enumerate}
		\item Both the ML and the MPL estimators are $\sqrt{n}$-consistent in the so called \textit{low temperature regime} (high values of the parameter $\beta$). This was first established in \cite{chspin} for the $p=2$ case, and later extended to tensor Ising models ($p > 2$) in \cite{jaesung1}. More precisely, it is shown in \cite {jaesung1} and \cite{jaesungcw} that there exists $\beta^*(p) > 0$, such that for all $\beta > \beta^*(p)$, both $\sqrt{n}(\bmp-\beta)$ and $\sqrt{n}(\bml - \beta)$ are tight, where $\bml$ is the Maximum Likelihood (ML) estimator of $\beta$. Further, consistent estimation (and consistent testing) is impossible in the regime $[0,\beta^*(p))$. The initial few values of the threshold $\beta^*(p)$ are given by $\beta^*(2)=0.5, \beta^*(3) = 0.672$ and $\beta^*(4)=0.689$. 
		
		\item Above the estimation threshold $\beta^*(p)$, both $\sqrt{n}(\bmp-\beta)$ and $\sqrt{n}(\bml - \beta)$ converge weakly to the same normal distribution \cite{jaesungcw}. This asymptotic normality can in fact be used to construct confidence intervals with asymptotically valid coverage probabilities for the parameter $\beta$ in presence of an external magnetic field term in the model \ref{cwmodel} (see Section 5 in \cite{jaesungcw}).
		
		\item As a consequence of the last point, both the ML and MPL estimates have the same asymptotic variance everywhere above the threshold. In fact, this asymptotic variance equals the limiting inverse Fisher information of the model, so both the estimates saturate the Cramer-Rao information lower bound of the model in this regime \cite{jaesungcw}.
	\end{enumerate} 
	
	What happens if we wish to perform a hypothesis testing of the natural parameter using these two estimators? In particular, which of the two tests would require a smaller sample size for achieving significance at a given level? The correct way to address this issue is to use the notion of \textit{Bahadur efficiency}.
	In this paper, we demonstrate that the MPL estimator is in fact as Bahadur efficient as the ML estimator everywhere in the classical $2$-spin $(p=2)$ model, and for \textit{most} values of the null and the alternative parameters in higher-order ($p\ge 3$) tensor models. When $p\ge 3$, a two-layer phase transition phenomenon is observed, depending on the magnitudes of the null and the alternative parameters. To elaborate, in this case, the only regime where the MPL estimator is less Bahadur efficient than the ML estimator, is a very small window of variation of the null parameter between $\beta^*(p)$ and $\log 2$, and that too, for large values of the alternative parameter only. Moreover, this small window $(\beta^*(p),\log 2)$ shrinks to the empty set as $p \rightarrow \infty$. This shows that the optimal sample size requirements for the tests based on the ML and the MPL estimators to achieve significance, are identical in the $2$-spin case and same over most of the parameter space in the higher-order tensor case.
	
	\begin{figure*}[h]\vspace{-0.12in}
		\centering
		\begin{minipage}[l]{1.0\textwidth}
			\centering
			\includegraphics[height=3.8in,width=6in]
			{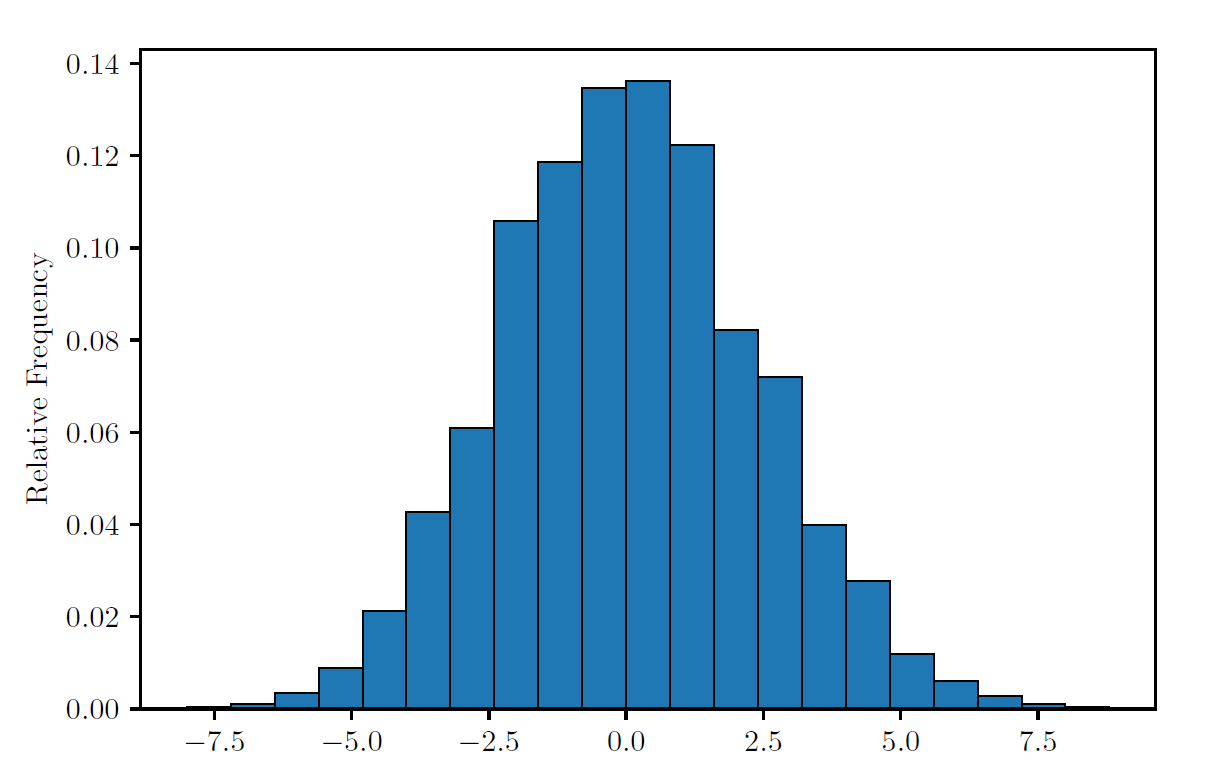}\\
		\end{minipage}
		\caption{\small{Histogram of $\sqrt n (\hat{\beta}_n(\bm X)-\beta)$, where $\hat{\beta}_n(\bm X)$ is the MPL estimator in the 4-tensor Curie-Weiss model at $\beta= 0.75 > \beta^*(4) \approx 0.689$ (above the estimation threshold), $n=20,000$ \cite{jaesung1}.}}
		\label{figpr}
	\end{figure*}
	In his seminal paper \cite{bahadur}, Bahadur introduced the concept of \textit{slope} of a test statistic to calculate the minimum sample size required to ensure its significance at a given level. The setting considered in \cite{bahadur} involved i.i.d. samples coming from a certain parametric family, and the goal was to detect the minimum sample size $N(\delta)$ required, so that a test $T_n$ (function of the samples) becomes (and remains) significant at level $\delta$ for all $n \ge N(\delta)$, i.e. the $p$-value corresponding to $T_n$ 
	becomes (and remains) bounded by $\delta$ for all $n \ge N(\delta)$. If one considers testing a simple null hypothesis $H_0: \theta = \theta_0$, then the above discussion may be quantified by defining:
	$$N(\delta) := \inf \left\{N \ge 1: \sup_{n \ge N} L_n \le \delta\right\}~,$$
	where $L_n:= 1-F_{T_n,\theta_0}(T_n)$ and $F_{T_n,\theta_0}$ is the cumulative distribution function of $T_n$ under $\p_{\theta_0}$. The $p$-value $L_n$ typically converges to $0$ exponentially fast with probability $1$ under alternatives $\p_{\theta}$ for $\theta > \theta_0$, and this rate is often an indication of the asymptotic efficiency of $T_n$ against $\theta$ \cite{anderson, bahadur2, bahadur2p5, bahadur3, bahadur4}. In particular, if we have the following $\p_\theta$-almost surely:
	\begin{equation}\label{pvalldp}
		\frac{1}{n} \log L_n \rightarrow -\frac{1}{2} c(\theta)\quad\textrm{as}~n \rightarrow \infty~,
	\end{equation} 
	then one can easily verify that (see Proposition 8 in \cite{bahadur}) $N(\delta) \sim -2\log(\delta)/c(\theta)$ as $\delta \rightarrow 0$. $c(\theta)$ is called the Bahadur slope of $T_n$ at $\theta$. However, as mentioned in \cite{bahadur}, it is in general a non-trivial problem to determine the existence of the Bahadur slope in \eqref{pvalldp}, and to evaluate it. This issue is addressed in two steps in \cite{bahadur}, where it is shown that if $T_n$ satisfies the following two conditions:
	\begin{enumerate}
		\item\label{i1} For every alternative $\theta$, $n^{-1/2} T_n \rightarrow b(\theta)$ as $n \rightarrow \infty$ under $\p_\theta$ with probability $1$, for some parametric function $b$ defined on the alternative space,
		
		\item\label{i2} $n^{-1}\log[1-F_{T_n,\theta_0}(n^{1/2} t)] \rightarrow - f(t)$ as $n \rightarrow \infty$ for every $t>0$ in an open interval which includes each value of $b$, where $f$ is a continuous function on the interval, with $0<f<\infty$,
	\end{enumerate}
	then the Bahadur slope exists for every alternative $\theta$, and is given by $2f(b(\theta))$ (see \cite{bahadur}). In this context, let us mention that if the convergence \eqref{pvalldp} holds in probability, then $c(\theta)$ is called the weak Bahadur slope of $T_n$ (see \cite{small_sample}). Finally, if we have two competing estimators $T_{n,1}$ and $T_{n,2}$ estimating the same parameter $\theta$, then the Bahadur asymptotic relative efficiency (ARE) is given by the ratio of their Bahadur slopes (see \cite{small_sample}):
	$$\mathrm{eff}(T_{n,1},T_{n,2}; \theta) = \frac {c_1(\theta)}{c_2(\theta)}~.$$ 
	
	The Bahadur ARE is a well-known tool in the literature for comparing the performance of two estimators in a wide variety of contexts. Gyorfi et al. \cite{gyorfi} addressed the problem of comparing the efficiencies of information-divergence-type statistics for testing the goodness of fit. They claim that the Pitman approach is too weak to detect sufficiently sharply the differences in efficiency of these statistics, and instead, focussed their attention on the Bahadur efficiency. Harremo\"{e}s and Vajda \cite{harremoes} show that in the problem of testing the uniformity of a distribution, the  information divergence  statistic is  more  efficient  in  the Bahadur  sense  than  any  power  divergence  statistic of order $\alpha > 1$. The same authors show in their paper \cite{havj} that any two R\'enyi entropies of different orders $\in (0,1]$ are equally Bahadur efficient. Huang \cite{huang} uses Bahadur efficiency as a measure of performance in the small sample universal hypothesis testing problem, and mentions that in the large sample problem where the number of possible outcomes is at most of the order of the number of samples, the connection between the error exponent and Bahadur efficiency has been studied in \cite{quine,harremoes}. Keziou and Regault \cite{kez} compare the performances of independence tests derived by means of dependence thresholding in a semiparametric context, in terms of the Bahadur ARE. Applications of the Bahadur efficiency in the contexts of nonparametric tests for independence and separate hypothesis testing can also be found in \cite{begh,ber} and \cite{rukhin}, respectively.

\begin{flushleft}
\begin{figure}
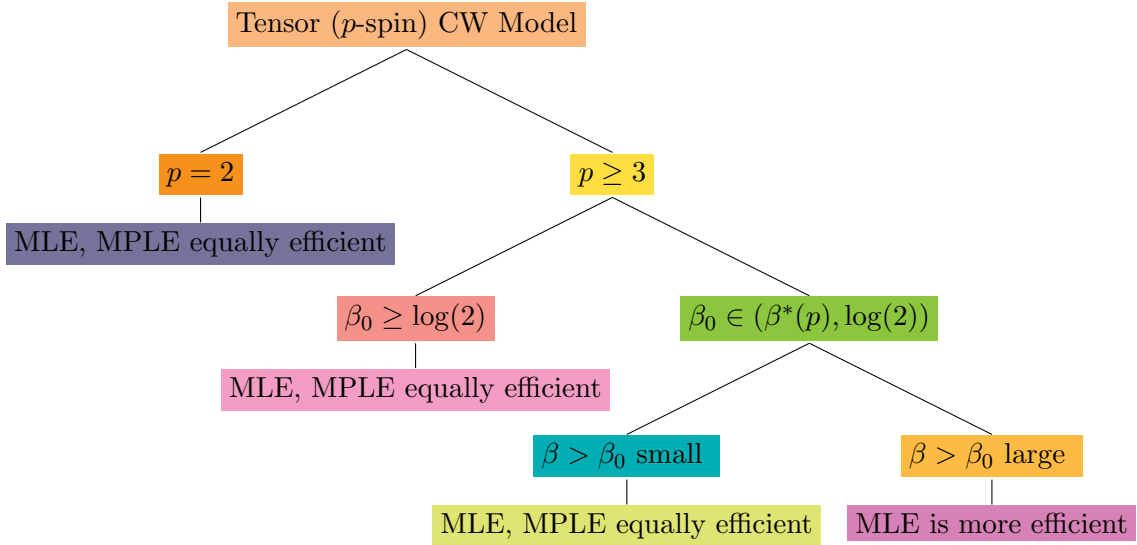

\Tree[.\colorbox{Apricot}{Tensor ($p$-spin) CW Model} [.\colorbox{BurntOrange}{$p=2$} \colorbox{CadetBlue}{MLE, MPLE  equally efficient}
               ]  !\qsetw{0.2cm}
          [.\colorbox{Goldenrod}{$p\geq3$} [.\colorbox{Salmon}{$\beta_0\geq\log(2)$} \colorbox{Lavender}{MLE, MPLE equally efficient} ] !\qsetw{0.1cm}
                [.\colorbox{LimeGreen}{$\beta_0\in(\beta^*(p),\log(2))$} [.\colorbox{TealBlue}{$\beta>\beta_0\text{ small }$} \colorbox{GreenYellow}{MLE, MPLE equally efficient} ]  
                [.\colorbox{Dandelion}{$\beta>\beta_0\text{ large }$} \colorbox{Thistle}{MLE is more efficient} ]]]]
\caption{Behavior of MLE and MPLE for different values of $p$}
\label{fig:mlempletree}
\end{figure}

\end{flushleft}



	
	
	In this paper, we compare the estimators $\bmp$ and $\bml$ in the tensor Curie-Weiss model \eqref{cwmodel} in terms of the Bahadur ARE. This requires  deriving the weak Bahadur slopes of both $\bmp$ and $\bml$ in the model \eqref{cwmodel}, which will in turn, enable one to compute the Bahadur ARE of either of these two estimators against some other reference estimator. Similar results have been derived in \cite{comets_lattice} in the context of Markov random fields on lattices, and in \cite{asymcomp} in the context of $d$-dimensional nearest neighbor isotropic Ising models, but to the best of our knowledge, this is the first such work on tensor Curie-Weiss (and Erd\H{o}s-R\'enyi Ising) models. Our basic tools are some recent results on large deviation of the \textit{average magnetization} $\os := \frac{1}{n} \sum_{i=1}^n X_i$ in the Curie-Weiss model, established in \cite{Lowe} and \cite{pspinldp}. Throughout the rest of the paper, we will view the entries $X_1,\ldots,X_n$ of the tuple $\bm X \in \{-1,1\}^n$ as \textit{dependent samples}, and refer to the length $n$ of $\bm X$ as the \textit{sample size} (although technically speaking, we have just one multivariate sample $\bm X$ from the model \eqref{cwmodel}). One of our most interesting findings is that the relative performance of the ML and the MPL estimators based on the notion of Bahadur efficiency depends crucially on whether the model is strictly tensor or not. In other words, in the usual $2$-spin ($p=2$) Curie-Weiss model, the two estimators are indistinguishable from the Bahadur ARE perspective. This is also true if the model is strictly tensor (i.e. $p\ge 3$), but for all values of the null parameter $\beta_0 \ge \log 2$. However, in this case, the MPL estimator is strictly less efficient than the ML estimator for each value of the null parameter $\beta_0$ in $\mathcal{R} := (\beta^*(p),\log 2)$, but that too, only if the alternative parameter $\beta$ is sufficiently large. This loss of Bahadur efficiency for the MPL estimator near threshold in the tensor Curie-Weiss models, can be attributed to its functional form, which derives false signal from a regime where the average magnetization $\os$ is very close to $0$. 
	
	Our results are a bit more universal, in the sense that they extend beyond the Curie-Weiss model \eqref{cwmodel}. They hold verbatim in the tensor \ern~ Ising model \eqref{ermodel} too, which is an exponential family, with sufficient statistic given by a tensor form, the tensor being the adjacency of a directed \ern~hypergraph with loops. Interestingly, we can even allow for slight sparsities in the underlying \ern~hypergraph. We believe that the same results (and techniques) will also extend to Ising models on dense stochastic block model hypergraphs, and leave it open for future research.
	
	The rest of the paper is organized as follows. In Section \ref{sec2}, we compare the ML and the MPL estimators in terms of their Bahadur ARE. The corresponding Bahadur slope and optimal sample size calculations necessary for the tests based on these two estimators are also provided in Section \ref{sec2}. Section \ref{sec:ermodel} shows that these results for the tensor Curie-Weiss model possess a somewhat universality property, i.e. they are true even for the tensor \ern~Ising model. In Section \ref{nr}, we provide numerical illustrations of our theoretical findings in various settings. In Section \ref{mainproofr}, we prove the main results in Sections \ref{sec2} and \ref{sec:ermodel}. In Section \ref{discussion}, we summarize some main and interesting aspects of our results, and talk about possible directions for future research in this area. Finally, proofs of some technical results needed for showing the main theorems are given in the appendix.
	
	\subsection{Our Contributions}
	Before delving in the main results, let us summarize our contributions once again in a pointwise fashion. We enlist the results we derived below.
	\begin{enumerate}
		\item In the classical $2$-spin $(p=2)$ Curie-Weiss/\ern~ model, the ML and MPL estimators are equally Bahadur efficient.
		
		\item In the higher-order tensor $(p\ge 3)$ Curie-Weiss/\ern~ model, the MPL estimator is equally Bahadur efficient as the ML estimator, if the null parameter is greater than or equal to $\log 2$. Even if the null parameter lies strictly between the model threshold and $\log 2$, they are equally efficient unless the alternative parameter is very large. In the last case, the MPL estimator is less Bahadur efficient than the ML estimator.
		
		\item For any arbitrary value of the null parameter above the threshold and any arbitrary value of the alternative parameter, a necessary and sufficient condition for the tests based on the ML and MPL estimators to be equally Bahadur efficient, is provided. 
		
		\item The exact Bahadur slopes of the ML and MPL estimators, and the optimal sample sizes required by the tests based on these estimators to achieve significance, are derived for both the tensor Curie-Weiss and the \ern~Ising model.
		
		\item Even if the the null parameter lies in the $(\beta^*(p),\log 2)$ window and the alternative parameter is very large, the MPL estimator is always at least a positive fraction as efficient as the ML estimator, an observation similar to Hodges and Lehmann's (1956) remarkable result that the Pitman ARE of Wilcoxon's test with respect to Student's T-test, under location alternatives, never falls below $0.864$, despite the former being non-parametric and exactly distribution-free for all sample sizes. The same conclusion holds for fixed values of the alternative above threshold, provided the null parameter is bounded away from the threshold.
		
		\item The Bahadur slopes and the optimal sample size requirements for both the ML and MPL estimators remain unchanged, if one moves from the tensor Curie-Weiss model to the tensor Erd\H{o}s-R\'enyi model. This indicates a possible universality of our results to Ising models on dense, random block hypergraphs.
		
	\end{enumerate} 
	
	Although the MPL estimator shows almost all the desirable efficiency properties of the ML estimator, in order to provide a completely unbiased and clear picture, we must mention two aspects in which the latter beats the former:
	
	\begin{enumerate}
		\item For $p\ge 3$ and every fixed value of the alternative parameter, the Bahadur ARE of the MPL estimator with respect to the ML estimator approaches $0$ as the null parameter approaches the threshold from the right. This is because, the expression for asymptotic optimal sample size \eqref{asmsamp} for significance of the test based on the MPL estimator approaches $\infty$, but the corresponding asymptotic optimal sample size for the test based on the ML estimator remains bounded.
		
		\item In the $p\ge 3$ case, as long as the null parameter lies strictly between the threshold and $\log 2$, the asymptotic optimal sample size required by the MPL estimator stabilizes at a fixed value after the alternative parameter exceeds a certain finite value, unlike that of the ML estimator. This is an undesirable property of the MPL estimator, since the sample size requirement does not decrease with increase in the separation between the null and alternative parameters, above a certain limit. However, as already mentioned above, the asymptotic optimal sample size for the MPL estimator exceeds that of the ML estimator by a bounded fraction only, as $\beta \rightarrow \infty$.
	\end{enumerate}
	
	\section{Theoretical Results for the Tensor Curie-Weiss Model}\label{sec2}
	
	In this section we compare the MPL and the ML estimates in the tensor Curie-Weiss model \eqref{cwmodel} in terms of their Bahadur ARE. The ML estimator $\bml$ does not have an explicit form, but it is shown in \cite{jaesung1} that the MPL estimator is given by:
	\[   
	\bmp = 
	\begin{cases}
		p^{-1} \os^{1-p} \tanh^{-1}(\os) &\quad\text{if} ~\os \ne 0,\\
		0 &\quad\text{if} ~\os = 0.\\
	\end{cases}
	\]
	Furthermore, it is shown in \cite{jaesungcw} and \cite{jaesung1} that both the ML and MPL estimators have the same asymptotic normal distribution:
	$$\sqrt{N}(\hat{\beta} - \beta)  \xrightarrow{D} N\left(0, -\frac{H_{\beta,p}''(m_*(\beta,p))}{p^2m_*(\beta,p)^{2p-2}}\right)~,$$ 
	for all $\beta > \beta^*(p)$, where $\hat{\beta}$ is either $\bml$ or $\bmp$,
	\begin{equation}\label{hdef}
		H_{\beta,p}(x) := \beta x^p - \frac{1}{2}\left\{(1+x)\log(1+x) + (1-x) \log(1-x) \right\} \quad \textrm{for}~ x \in [-1,1]~,
	\end{equation}
	$m_*(\beta,p)$ is the unique positive global maximizer of $H_{\beta,p}$, and 
	$$\beta^*(p) := \sup\left\{\beta \ge 0: \sup_{x\in [-1,1]} H_{\beta,p}(x) = 0\right\}.$$ A few initial values of the threshold $\beta^*(p)$ are $\beta^*(2)= 0.5, \beta^*(3)\approx 0.672$ and $\beta^*(4) \approx 0.689$. The exact value of $\beta^*(p)$ is in general inexplicit, but $\beta^*(p) \uparrow \log 2$ as $p \rightarrow \infty$ (see Lemma A.1 in \cite{jaesung1}). 
	
	In this paper, we will consider testing the hypothesis $$H_0: \beta = \beta_0\quad \quad \textrm{vs}\quad \quad H_1: \beta > \beta_0$$ for some known $\beta_0 > \beta^*(p)$. The most powerful test for this hypothesis in model \eqref{cwmodel} is based on the sufficient statistic $\os$, and its asymptotic power is derived in \cite{BM16}. Clearly, one can think of using the statistic $T_n:= \sqrt{n}(\hat{\beta} -\beta_0)$ for testing the above hypotheses, where $\hat{\beta}$ is either $\bml$ or $\bmp$, and large values of $T_n$ will denote significance.
	
	We now state the main result in this paper about the Bahadur slopes of the tests based on the MPL and ML estimators, and the minimum sample size required to ensure their significance. Towards this, we define a function $\eta_p: [-1,1]\mapsto \mathbb{R}$ as:
	\[   
	\eta_p(t) = 
	\begin{cases}
		p^{-1} t^{1-p} \tanh^{-1}(t) &\quad\text{if} ~t \ne 0,\\
		0 &\quad\text{if} ~t = 0.\\
	\end{cases}
	\]
	
	\begin{theorem}\label{main_mpl}
		The Bahadur slopes of $\bmp$ and $\bml$ for the model \eqref{cwmodel} at an alternative $\beta$ are respectively given by:
		$$c_{\bmp}(\beta_0,\beta,p) = 2\left(\sup_{x\in [-1,1]} H_{\beta_0,p}(x) - \sup_{x \in \eta_p^{-1} ((\beta,\infty))} H_{\beta_0,p}(x)\right)~,$$
		$$c_{\bml}(\beta_0,\beta,p) = 2\left(\sup_{x\in [-1,1]} H_{\beta_0,p}(x) - \sup_{x ~>~ m_*(\beta,p)} H_{\beta_0,p}(x)\right)$$
		Consequently, the minimum sample sizes required, so that the tests $\sqrt{n}(\bmp-\beta_0)$ and $\sqrt{n}(\bml - \beta_0)$ become (and remain) significant at level $\delta \rightarrow 0$, are respectively given by:
		$$N_{\bmp}(\beta_0,\beta,\delta,p) ~\sim~ \frac{\log(\delta)}{\sup_{x \in \eta_p^{-1} ((\beta,\infty))} H_{\beta_0,p}(x)-\sup_{x\in [-1,1]} H_{\beta_0,p}(x)}~,$$
		$$N_{\bml}(\beta_0,\beta,\delta,p)~ \sim ~\frac{\log(\delta)}{\sup_{x ~>~ m_*(\beta,p)} H_{\beta_0,p}(x) - \sup_{x\in [-1,1]} H_{\beta_0,p}(x) }~.$$
	\end{theorem}

	
	
	
	

	Theorem \ref{main_mpl} is proved in Section \ref{mainmplproof7}. Let us introduce the following notation, which will be used throughout the rest of the paper:
	\begin{equation}\label{asmsamp}
		N_{\bmp}^*(\beta_0,\beta,\delta,p)= -\frac{2\log(\delta)}{c_{\bmp}(\beta_0,\beta,p)}\quad\text{and}\quad N_{\bml}^*(\beta_0,\beta,\delta,p)= -\frac{2\log(\delta)}{c_{\bml}(\beta_0,\beta,p)}~.
	\end{equation}
	Note that $N_{\bmp}(\beta_0,\beta,\delta,p) \sim N_{\bmp}^*(\beta_0,\beta,\delta,p)$ and $N_{\bml}(\beta_0,\beta,\delta,p) \sim N_{\bml}^*(\beta_0,\beta,\delta,p)$ as $\delta \rightarrow 0$, so the latter quantities can be called the \textit{asymptotic} optimal sample sizes.
	The natural question at this stage, is when do the expressions for the Bahadur slope and the asymptotic optimal sample size for the MPL and the ML estimators differ? It turns out that we have two different scenarios depending on whether $p=2$ or $p\ge 3$. 
	
	\subsection{The $p=2$ Case:}
	In this case, it turns out that the Bahadur slopes and the asymptotic optimal sample sizes for the MPL and the ML estimators agree, and consequently, they are equally Bahadur efficient.
	
	\begin{theorem}\label{p2case}
		Consider the model \eqref{cwmodel} for $p=2$. For every $\beta>\beta_0>\beta^*(2)$ and $\delta \in (0,1)$, we have $$c_{\bmp}(\beta_0,\beta,2) = c_{\bml}(\beta_0,\beta,2)\quad\textrm{and}\quad N_{\bmp}^*(\beta_0,\beta,\delta,2) = N_{\bml}^*(\beta_0,\beta,\delta,2)~.$$ Consequently, the Bahadur ARE $\mathrm{eff}(\bml,\bmp;\beta_0,\beta) = 1$.
	\end{theorem}
	Theorem \ref{p2case} is proved in Section \ref{p2cproof}. It says that in the classical $2$-spin Curie-Weiss model, one cannot distinguish the estimators $\bml$ and $\bmp$ based on even the Bahadur efficiency. 
	
	\subsection{The Strictly Tensor $p \ge 3$ Case:}
	
	All the interesting phenomena occur when the model \eqref{cwmodel} goes beyond the classical $2$-spin system to the higher-order tensor $(p\ge 3)$ system. In this case, a two layer phase transition is observed with respect to both the null and alternative parameters. To be precise, there exists a small window around the estimation threshold $\beta^*(p)$, such that for all values of the null parameter in this window, $\bmp$ is strictly less Bahadur efficient than $\bml$ provided the alternative parameter $\beta$ is greater than a second threshold (which is different from $\beta^*(p)$).
	
	\begin{theorem}\label{p3case}
		Consider the model \eqref{cwmodel} for $p\ge 3$. Two different situations arise depending on whether $\beta_0 \ge \log 2$ or $\beta_0 \in (\beta^*(p),\log 2)$.
		
		\begin{enumerate}
			\item For every $\beta>\beta_0\ge \log 2$ and $\delta \in (0,1)$, we have $$c_{\bmp}(\beta_0,\beta,p) = c_{\bml}(\beta_0,\beta,p)\quad\textrm{and}\quad N_{\bmp}^*(\beta_0,\beta,\delta,p) = N_{\bml}^*(\beta_0,\beta,\delta,p)~.$$ Consequently, the Bahadur ARE $\mathrm{eff}(\bml,\bmp;\beta_0,\beta) = 1$ in this regime.
			\vspace{0.05in}
			\item For every $\beta_0 \in (\beta^*(p),\log 2)$ and $\delta \in (0,1)$, the Bahadur slopes and asymptotic optimal sample sizes for $\bmp$ and $\bml$ do not agree, and $\mathrm{eff}(\bml,\bmp;\beta_0,\beta) > 1$, for all $\beta>\beta_0$ large enough.
		\end{enumerate}
	\end{theorem}
	
	Theorem \ref{p3case} is proved in Section \ref{p3casepr7}. Some more insight is obtained if one fixes the alternative $\beta > \beta^*(p)$, and looks at the behavior of the asymptotic optimal sample sizes and Bahadur ARE of the MPL and ML estimators by varying the null parameter in the small window near $\beta^*(p)$. 
	
	\begin{theorem}[The fixed alternative scenario]\label{beta1fixed}
		For $p\ge 3$ and fixed $\beta > \beta^*(p)$, the set of all $\beta_0 > \beta^*(p)$ satisfying $\mathrm{eff}(\bml,\bmp;\beta_0,\beta) > 1$ is given by:
		$$\left(\beta^*(p),\frac{I(m_*(\beta,p))}{m_*(\beta,p)^p}\right)$$ where $I(x) = \frac{1}{2}\left\{(1+x)\log(1+x) + (1-x) \log(1-x) \right\}$. Further, for every $p\ge 3$ and every fixed $\beta>\beta^*(p)$, we have:
		\begin{equation}\label{last}
			\lim_{\beta_0\rightarrow \beta^*(p)^+} N_{\bmp}^*(\beta_0,\beta,\delta,p) = \infty \quad\textrm{and}\quad\lim_{\beta_0\rightarrow \beta^*(p)^+} N_{\bml}^*(\beta_0,\beta,\delta,p) = \frac{\log(\delta)}{H_{\beta^*(p),p}(m_*(\beta,p))} < \infty.
		\end{equation}
		and hence,
		$$\lim_{\beta_0\rightarrow \beta^*(p)^+} \mathrm{eff}(\bml,\bmp;\beta_0,\beta) = \infty~.$$
	\end{theorem}
	
	Theorem \ref{beta1fixed} is proved in Section \ref{betafixedpr7}.
	
	\begin{remark}
		Theorem \ref{beta1fixed} says that as the null parameter approaches the estimation threshold, the ML estimator becomes infinitely more Bahadur efficient than the MPL estimator. The reason behind the discrepancy between the efficiencies of the ML and MPL estimators near the threshold, is the functional form of the latter. For $p\ge 3$, unlike the ML estimator, the MPL estimator takes very high values if the average magnetization $\os$ is close to $0$. This false signal coming from the average magnetization lying in a region very close to $0$, leads to an increase in the null probability of the MPL estimator exceeding the observed MPL estimate, thereby inflating its $p$-value. This inflation occurs only in a close neighborhood of the threshold, because for lower values of the parameter $\beta$, there is a higher probability that the average magnetization $\os$ is small. 
	\end{remark}
	
	\begin{remark}
		It follows from Lemma C.12 in \cite{jaesungcw} that 
		$$\sup_{x \in [-1,1]}H_{\beta_0,p}(x) = H_{\beta_0,p}(m_*(\beta_0,p)) = \Theta(\beta_0-\beta^*(p))~.$$ Since $\sup_{x \in \eta_p^{-1}((\beta,\infty))} H_{\beta_0,p}(x)$ eventually becomes $0$ as $\beta_0$ approaches $\beta^*(p)$ from the right, the rate at which $N_{\bmp}^*(\beta_0,\beta,\delta,p)$ approaches $\infty$ as $\beta_0 \rightarrow \beta^*(p)^+$ for $p \ge 3$, is determined just by the $\sup_{x \in [-1,1]}H_{\beta_0,p}(x)$ term in the denominator of the formula for $N_{\bmp}^*(\beta_0,\beta,\delta,p)$, and is given by $(\beta_0-\beta^*(p))^{-1}$.
	\end{remark}
	
	Another undesirable property of $\bmp$ is that for every fixed value of the null parameter in the window $(\beta^*(p),\log 2)$, its asymptotic optimal sample size does not decrease further when the alternative parameter exceeds a certain value. So, no matter how large the separation between the null and the alternative are, one does not have any concession in the asymptotic sample size requirement after a certain value of the separation. This is formalized in the theorem below:
	
	\begin{theorem}\label{stablesamp}
		For $p\ge 3$ and fixed $\beta_0 \in (\beta^*(p),\log 2)$, there exists $\underline{\beta} > 0$ such that $$N_{\bmp}^*(\beta_0,\beta,\delta,p) = -\frac{\log(\delta)}{H_{\beta_0,p}(m_*(\beta_0,p))} > 0\quad\text{for all} ~\beta > \underline{\beta}.$$
	\end{theorem}
	
	Theorem \ref{stablesamp} is proved in Section \ref{th5finpr}. It says that the asymptotic optimal sample size requirement for the MPL estimator stabilizes at a certain value once the separation between the null and the alternative is large enough, provided the null lies in the window $(\beta^*(p),\log 2)$. In this case, note that for $\beta$ large enough, we have:
	$$N_{\bml}^*(\beta_0,\beta,\delta,p) = -\frac{\log(\delta)}{H_{\beta_0,p}(m_*(\beta_0,p)) - H_{\beta_0,p}(m_*(\beta,p))} < N_{\bmp}^*(\beta_0,\beta,\delta,p).$$
	It also shows that unlike the $\bmp$, the asymptotic optimal sample size for $\bml$ is strictly decreasing with the alternative $\beta$, for all sufficiently large values of $\beta$. This is precisely due to the presence of the $H_{\beta_0,p}(m_*(\beta,p))$ term in the denominator of $N_{\bml}^*(\beta_0,\beta,\delta,p)$, which is strictly decreasing in $\beta$ for all $\beta$ large enough.
	
	\begin{remark}\label{Tbounded}
		It follows from Theorem \ref{main_mpl} that for every fixed $\beta_0 > \beta^*(p)$,
		$$\lim_{\beta\rightarrow \infty} N_{\bml}^*(\beta_0,\beta,\delta,p) = \frac{\log(\delta)}{\beta_0-\log 2 - H_{\beta_0,p}(m_*(\beta_0,p))}~.$$ On the other hand, Theorem \ref{stablesamp} implies that as long as $\beta_0 \in (\beta^*(p),\log 2)$, we have:
		$$\lim_{\beta\rightarrow \infty} N_{\bmp}^*(\beta_0,\beta,\delta,p) = \frac{\log(\delta)}{- H_{\beta_0,p}(m_*(\beta_0,p))}~.$$ Hence, for every $\beta_0 \in (\beta^*(p),\log 2)$,
		$$\lim_{\beta\rightarrow \infty} \mathrm{eff}(\bmp,\bml;\beta_0,\beta) = \frac{H_{\beta_0,p}(m_*(\beta_0,p))}{H_{\beta_0,p}(m_*(\beta_0,p)) +\log 2-\beta_0} > 0~.$$ Hence, the Bahadur ARE of the MPL estimator with respect to the ML estimator is bounded away from $0$ even when the null parameter is in the interval $(\beta^*(p),\log 2)$ and the alternative parameter is large. This is analogous to the remarkable result of Hodges and Lehmann (1956), which shows that the Pitman ARE of Wilcoxon's test with respect to Student's T-test, under location alternatives, never falls below $0.864$, despite the former being non-parametric and exactly distribution-free for all sample sizes (see \cite{bnb}).
	\end{remark}

	\section{Theoretical Results for the Hypergraph Erd\H{o}s-R\'enyi Ising Model}\label{sec:ermodel}
	In this section, we are going to see that all our results for the tensor Curie-Weiss model hold verbatim for the more general class of hypergraph Erd\H{o}s-R\'enyi Ising model, where we can even introduce some sparsity. A classical survey on these and more general models of disordered ferromagnets in the physics literature can be found in \cite{kab1}. This model for $p=2$ was introduced and studied in \cite{bovierjsp}, and analyzed further in \cite{erm1}, where a quenched central limit theorem has been proved for the magnetization in the high temperature (small $\beta$) regime, allowing for some sparsity in the underlying random graph. Of course, we will be concerned with the general $p$ case in this section.
	
	Let $\bm A:= \{A_{i_1\ldots i_p}\}_{1\leq i_1,\ldots,i_p \leq n}$ be a collection of i.i.d. Bernoulli random variables with mean $\alpha_n$. Note that $\bm A$ can be viewed as the adjacency tensor of a directed Erd\H{o}s-R\'enyi hypergraph with loops. The $p$-tensor \er~ Ising model in this context, is a discrete exponential family on $\{-1,1\}^n$ with probability mass function given by:
	\begin{equation}\label{ermodel}
		\p_{\beta,p}^*(\bt) = \frac{\exp\{\beta H_n(\bt)\}}{2^n Z_n^*(\beta,p)}\quad(\text{for}~\bt \in \{-1,1\}^n)~,
	\end{equation}
	where $$H_n(\bt) := \alpha_n^{-1} n^{1-p} \sum_{1\leq i_1,\ldots,i_p \leq n} A_{i_1\ldots i_p} x_{i_1}\ldots x_{i_p}$$ denotes the Hamiltonian of the model, and $Z_n^*(\beta,p)$ is the normalizing constant. Note that we will use a $*$ superscript to denote probabilities and moments corresponding to the model \eqref{ermodel}. Below, we state the main result of this section:
	
	\begin{theorem}\label{ermain}
		The Bahadur slopes and asymptotic minimum sample sizes of $\bmp$ and $\bml$ for the model \eqref{ermodel} at an alternative $\beta$ are respectively equal to the Bahadur slopes and asymptotic minimum sample sizes of $\bmp$ and $\bml$ for the model \eqref{cwmodel} at $\beta$.
	\end{theorem}
	
	Consequently, all the results in Section \ref{sec2} for the tensor Curie-Weiss model also hold for the tensor \ern~Ising model. Theorem \ref{ermain} is proved in Section \ref{erproof7}. The main approach, in a nutshell, is an approximation of the Hamiltonian and the local fields of the tensor \ern~model by the corresponding quantities for the tensor Curie-Weiss model. In order to show that the Hamiltonian $H_n(\bt)$ of the hypergraph \ern~Ising model is very close to that of the tensor Curie-Weiss model, it is enough to establish a uniform (in $\bt$) concentration of $H_n(\bt)$ around its mean (with respect to the \ern~measure) $\e H_n(\bt)$. Define $$\gamma_n := 3(\alpha_n n^{p-1})^{-\frac{1}{2}}~.$$
	
	\begin{lemma}\label{hamiltcomp}
		Let $H_n$ denote the Hamiltonian of the $p$-tensor \ern~Ising model.
		Then,
		$$\p\left(\frac{1}{n}\sup_{\bm x \in \sa} |H_n(\bt) - \e H_n(\bt)| \leq 3\gamma_n~\textrm{for all but finitely many}~n\right) =1~.$$
	\end{lemma}
	
	Lemma \ref{hamiltcomp} says that as long as the \ern~hyperedge probability $\alpha_n \gg n^{1-p}$, the Hamiltonian $H_n(\bt)$ concentrates around its mean (when scaled by a factor of $1/n$) uniformly in $\bt$.   The proof of Lemma \ref{hamiltcomp} is given in Appendix \ref{ertech}. The following result is a corollary of the proof of Lemma \ref{hamiltcomp}, which says that not only are the Hamiltonians of the two models close, but their local fields $m_i^{(n)}(\bt)$ (defined below) are also close, uniformly over all $1\le i\le n$ and all $\bt \in \{-1,1\}^n$. It will be useful in deriving the Bahadur slope of the MPL estimator.
	
	\begin{corollary}\label{hcm1}
		For each $1\le i\le n$, define $m_i^{(n)}(\bm x) := \alpha_n^{-1} n^{1-p} \sum_{(i_2,\ldots,i_p)\in [n]^{p-1}} A_{ii_2\ldots i_p} x_{i_2}\ldots x_{i_p}$.
		Then, we have:
		$$\p\left(\sup_{1\le i\le n} ~\sup_{\bm x \in \sa} \left|m_i^{(n)}(\bm x) - \oss^{p-1}\right| \leq 3\gamma_n~\textrm{for all but finitely many}~n\right) =1~.$$
	\end{corollary}
	
	The proof of Corollary \ref{hcm1} is given in Appendix \ref{ertech}. We will henceforth assume the slightly stronger condition $\alpha_n = \Omega(n^{1-p} \log n)$, which in particular implies that $\gamma_n \ll 1$. Note that for $p=2$, this condition is satisfied if the Erd\H{o}s-R\'enyi graph is almost-surely connected. This technical condition is required for establishing the consistency of the MPL estimator (see Lemma \ref{c11}). See Section \ref{erproof7} for a complete proof of Theorem \ref{ermain}.

	\subsection{Sketch of Proofs of the Main Results}
	We now provide a brief sketch of some of the main ideas involved in the proofs of the main results stated in Sections \ref{sec2} and \ref{sec:ermodel}. The main result in Section \ref{sec2} is Theorem \ref{main_mpl}, and its proof starts by establishing the probability limit of $\hat{\beta} -\beta_0$ under the alternative, where $\hat{\beta}$ denotes either the ML or the MPL estimator (see Lemma \ref{cond1}). This verifies Condition \eqref{i1}. In order to verify Condition \eqref{i2}, we use a large deviation principle of $\overline{X}_n$ to pinpoint the asymptotic behavior of the $p$-values (see Lemma \ref{xbarldp}). The proofs of the other results in Section \ref{sec2} follow from Theorem \ref{main_mpl} and an analysis of the behavior of the function $H_{\beta_0,p}$ and its maximum value over certain regions. 
	
	The proof of Theorem \ref{ermain} relies on comparing the Erd\H{o}s-R\'enyi Ising model with the Curie-Weiss model, and showing that the two models are close in the exponential scale (see Lemma \ref{measureclose}). This follows from an exponential concentration of the Hamiltonian of the Erd\H{o}s-R\'enyi Ising model around its expected value, the latter being the Hamiltonian of the Curie-Weiss model. This keeps the asymptotics of the parameter estimates and the large deviation behavior of the sample mean unchanged in the Erd\H{o}s-R\'enyi Ising model.

	\section{Numerical Results}\label{nr}

	In this section, we provide a graphical presentation of the numerical values of the asymptotic optimal sample sizes for the tests based on the ML and MPL estimators in the model \eqref{cwmodel}, using the theoretical formula given in Theorem \ref{main_mpl}. Note that the Curie-Weiss model can be rewritten as
	\begin{equation*}
		\p_{\beta, h,p} (\overline{X}_N = m) =\frac{1}{2^n Z_n(\beta, p)}\binom{n}{n(1+m)/2} e^{n\beta m^p}
	\end{equation*}
	for $m\in \left\{-1,-1 + \frac{2}{n}, \cdots, 1- \frac{2}{n}, 1\right\}$, whence the partition function can be computed as
	\begin{equation*}
		Z_n(\beta,p) = \sum_{m\in \mathcal{M}} \frac{1}{2^n} \binom{n}{n(1+m)/2} e^{n \beta m^p}.
	\end{equation*}
	Consequently, the mass function can be computed exactly for all moderately large $n$ (say $n$ up to $1000$) and the Curie-Weiss model can be sampled directly from this mass function without the use of any MCMC approach.
	
\begin{figure}
		\centering
		\includegraphics[height=3in,width=4in]{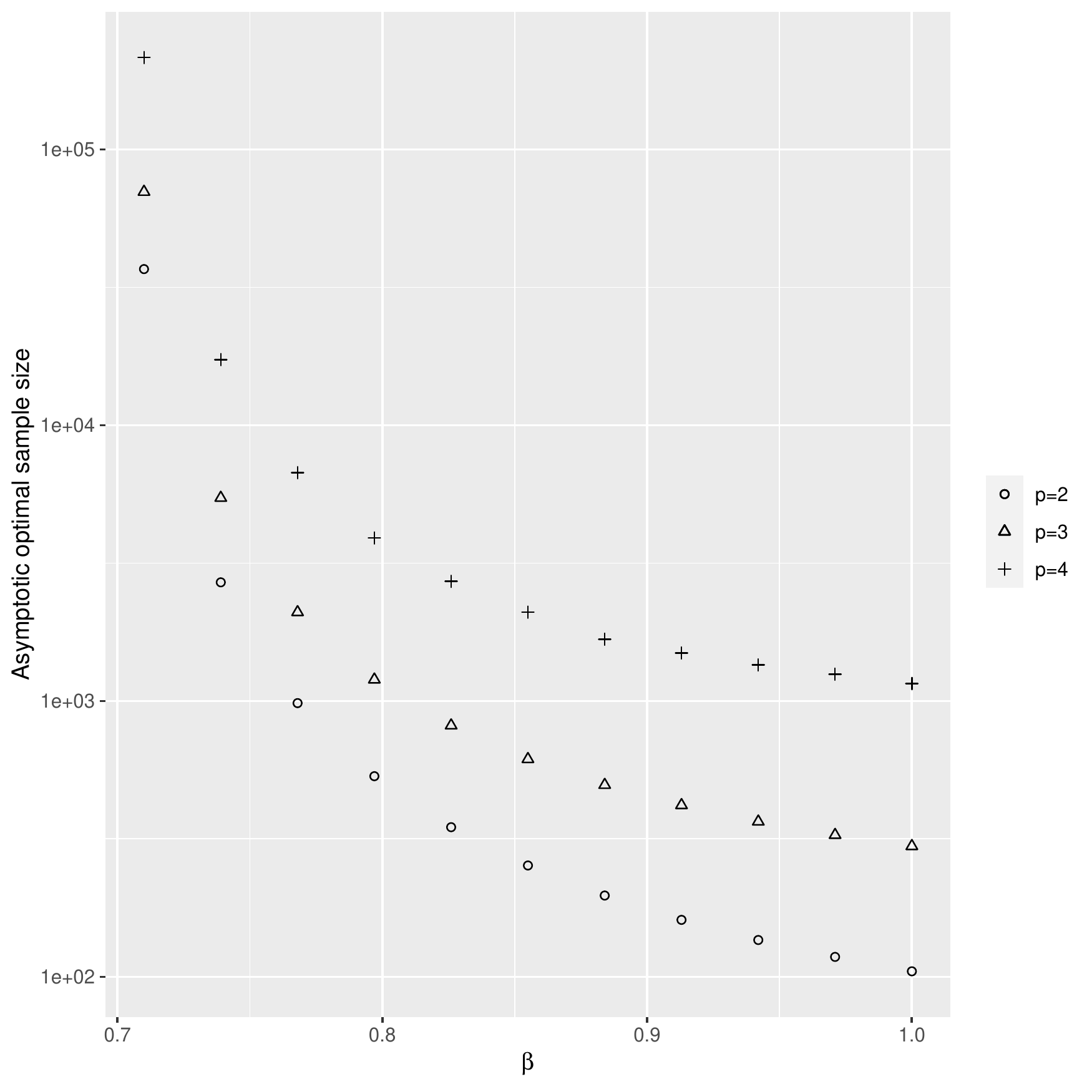}
		\caption{\small{Asymptotic optimal sample size for the tests based on MLE and MPLE with varying $\beta; ~\beta_0 = 0.7>\log 2,\delta=0.05, p\in \{2,3,4\}$ (with logarithmic vertical scale).}}
		\label{fig00}
	\end{figure}

	In Figures 3--7, we fix the level $\delta = 0.05$. In Figure \ref{fig00}, we fix $\beta_0 =0.7$, a value slightly larger than $\log 2$, and plot the asymptotic optimal sample size (for both the MPL and the ML tests, which must be same in this regime) for $p=2,3$ and $4$, across $\beta> \beta_0$. We see that the asymptotic optimal sample size decreases as the alternative parameter $\beta$ increases, which is expected, since the detection capability of the tests should increase as the alternative $\beta$ moves far apart from the null $\beta_0$. Another important observation is that with increase in the interaction complexity $p$, the asymptotic optimal sample size requirement also increases at every alternative. One possible explanation of this phenomenon is that with increase in $p$, the threshold $\beta^*(p)$ also increases, and hence, the null $\beta_0$ (which is fixed at 0.7) gets closer to $\beta^*(p)$, which causes a slight increase in the difficulty of the testing problem.
	
	\begin{figure}
		\centering
		
		\begin{minipage}{.5\textwidth}
			\centering
			
			\includegraphics[width=1\linewidth, height=0.28\textheight]{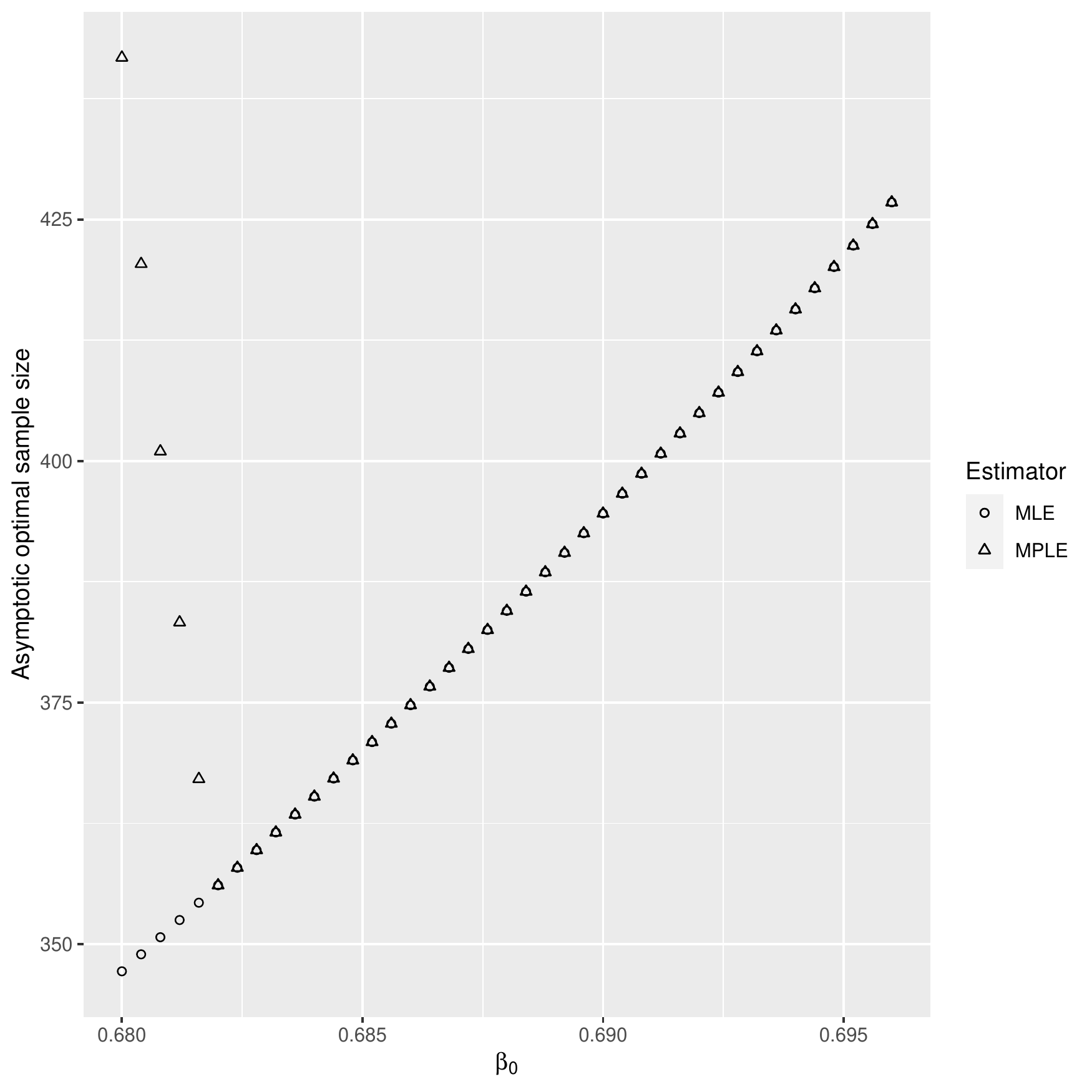}
			\caption{\small{Asymptotic optimal sample size for the tests based on MLE and MPLE with varying $\beta_0;~p=3,\beta=0.90,\delta=0.05$.}}
			\label{fig2}
		\end{minipage}%
		\begin{minipage}{0.5\textwidth}
			\centering
			\includegraphics[width=1\linewidth, height=0.28\textheight]{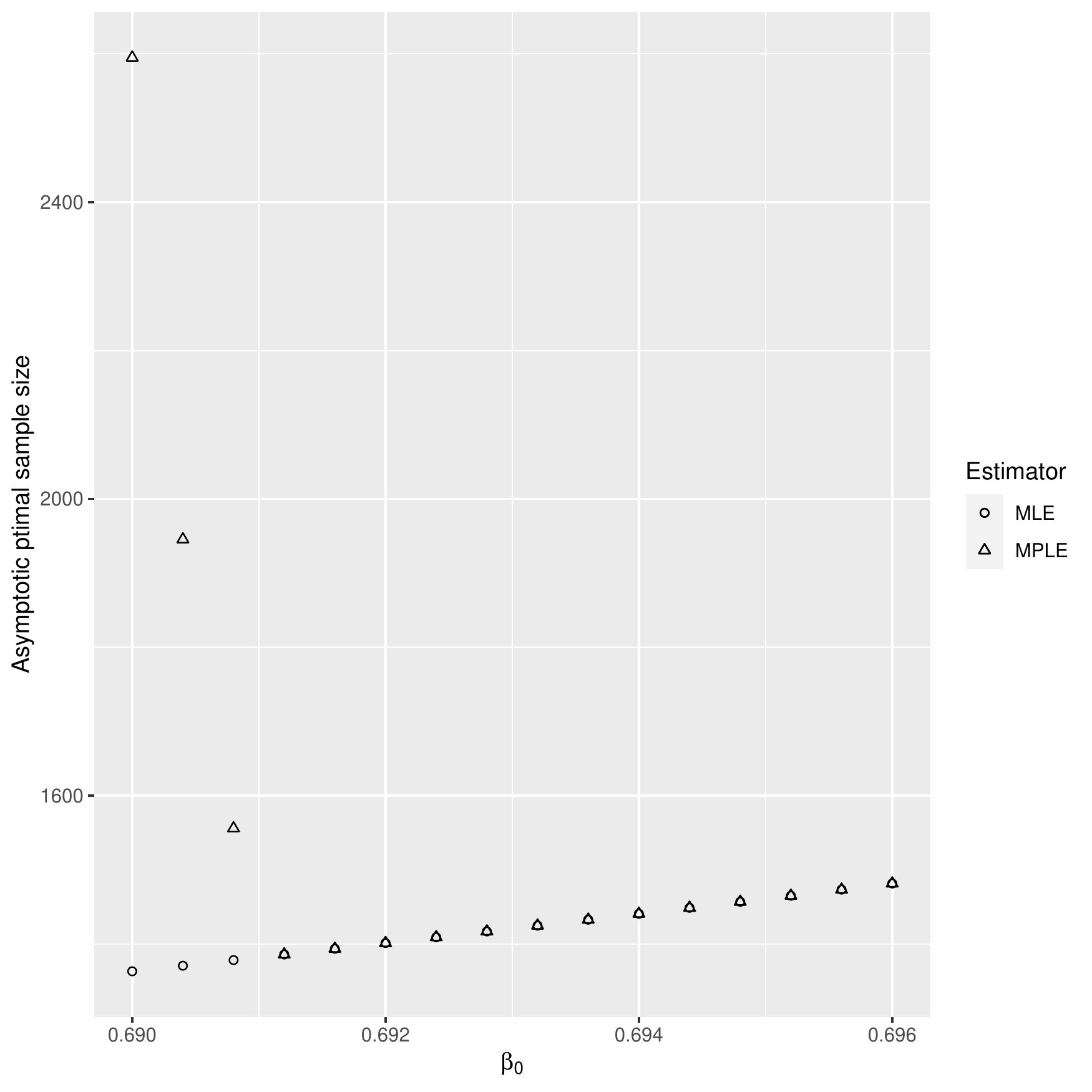}
			\caption{\small{Asymptotic optimal sample size for the tests based on MLE and MPLE with varying $\beta_0;~ p=4,\beta=0.90,\delta=0.05$.}}
			\label{fig3}
		\end{minipage}
		
	\end{figure}
	
\begin{figure}
		\centering
		
		\begin{minipage}{.5\textwidth}
			\centering
			\includegraphics[width=1\linewidth, height=0.28\textheight]{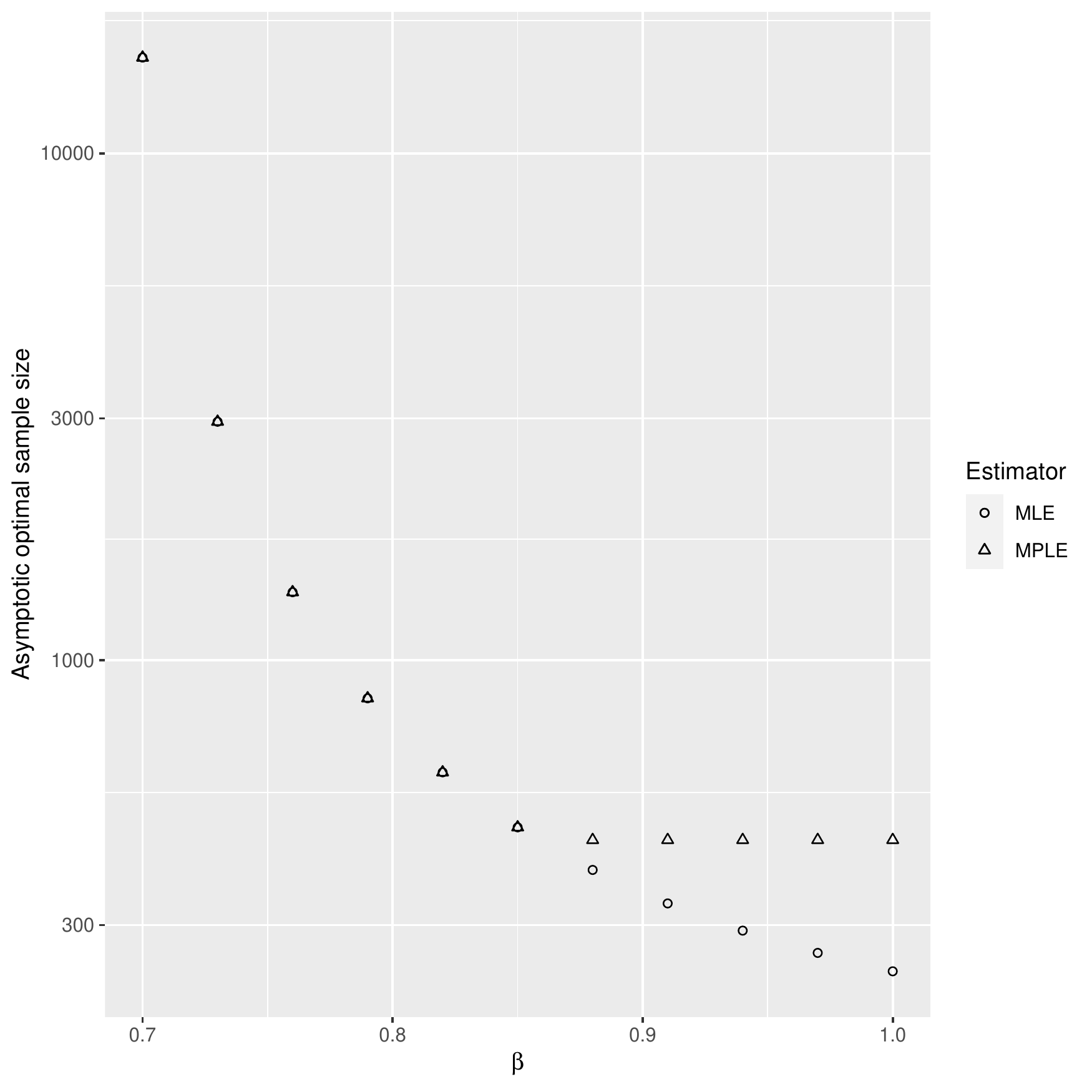}
			\caption{\small{Asymptotic optimal sample size for the tests based on MLE and MPLE with varying $\beta;~ p=3,\beta_0=0.68 < \log 2,\delta=0.05$ (with logarithmic vertical scale).}}
			\label{fig4}
		\end{minipage}%
		\begin{minipage}{0.5\textwidth}
			
			\centering
			\includegraphics[width=1\linewidth, height=0.28\textheight]{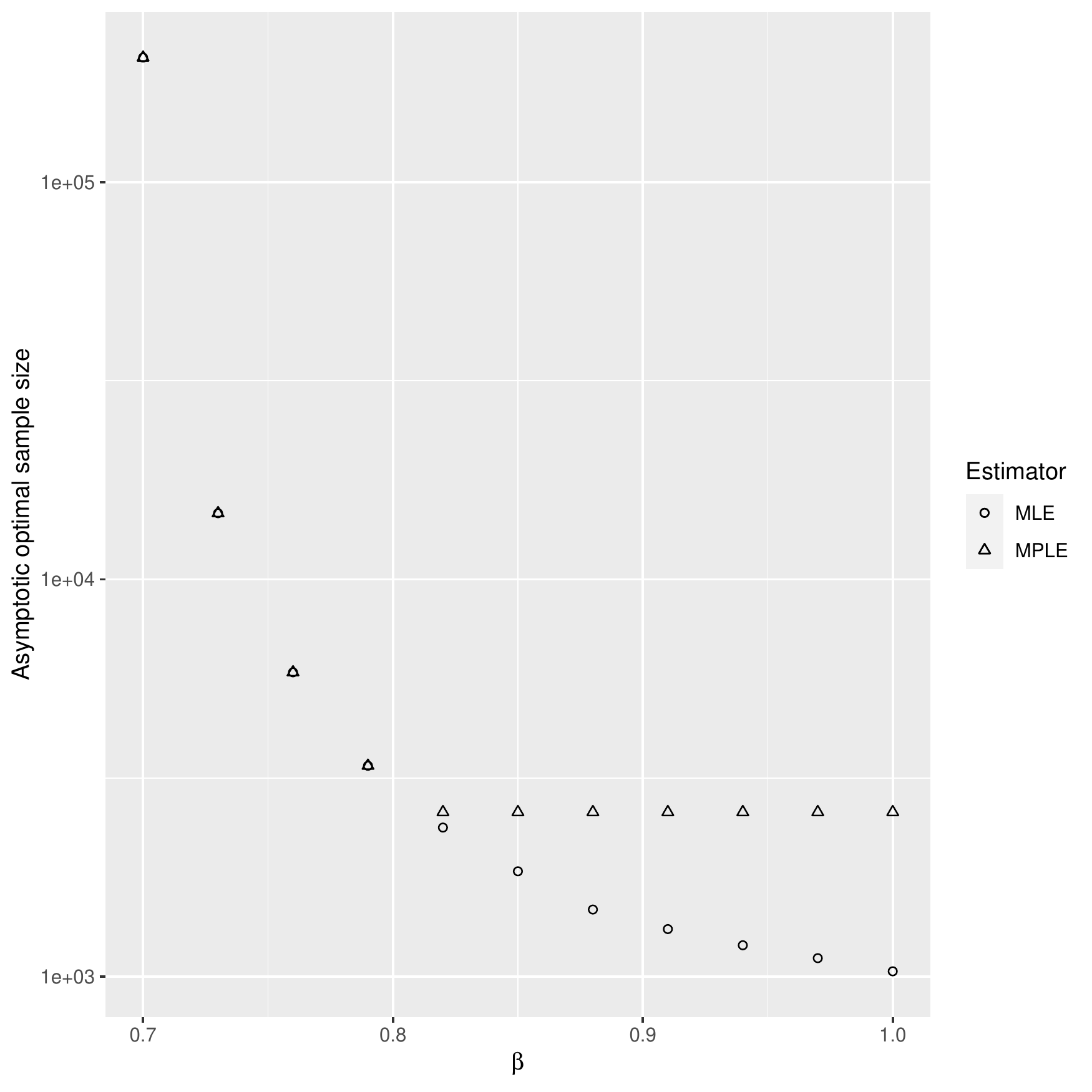}
			\caption{\small{Asymptotic optimal sample size for the tests based on MLE and MPLE with varying $\beta;~p=4,\beta_0=0.69 <\log 2,\delta=0.05$ (with logarithmic vertical scale).}}
			\label{fig5}
		\end{minipage}
		
	\end{figure}

	In Figures \ref{fig2} and \ref{fig3}, we fix the alternative $\beta=0.9$, and demonstrate graphically (for the cases $p=3$ and $4$ respectively), that the asymptotic optimal sample sizes for the MPL and the ML tests differ for all $\beta_0$ in a small right neighborhood of $\beta^*(p)$ below $\log 2$, and agree above that neighborhood. The figures also demonstrate that the asymptotic optimal sample size for the MPL test approaches $\infty$ as the null $\beta_0$ approaches the threshold $\beta^*(p)$. In Figures \ref{fig4} and \ref{fig5}, we fix the null $\beta_0$ to values slightly smaller than $\log 2$, and demonstrate (for the cases $p=3$ and $4$ respectively) that although the asymptotic optimal sample sizes for the MPL and the ML tests coincide for all small values of the alternative $\beta > \beta_0$, they disagree for all $\beta$ large enough. All these results reflect the contents of Theorems \ref{p3case} and \ref{beta1fixed}. 

	
	
\begin{figure}
		\centering
		
		\begin{minipage}{.5\textwidth}
			\centering
			\centering
			\includegraphics[width=1\linewidth, height=0.28\textheight]
			{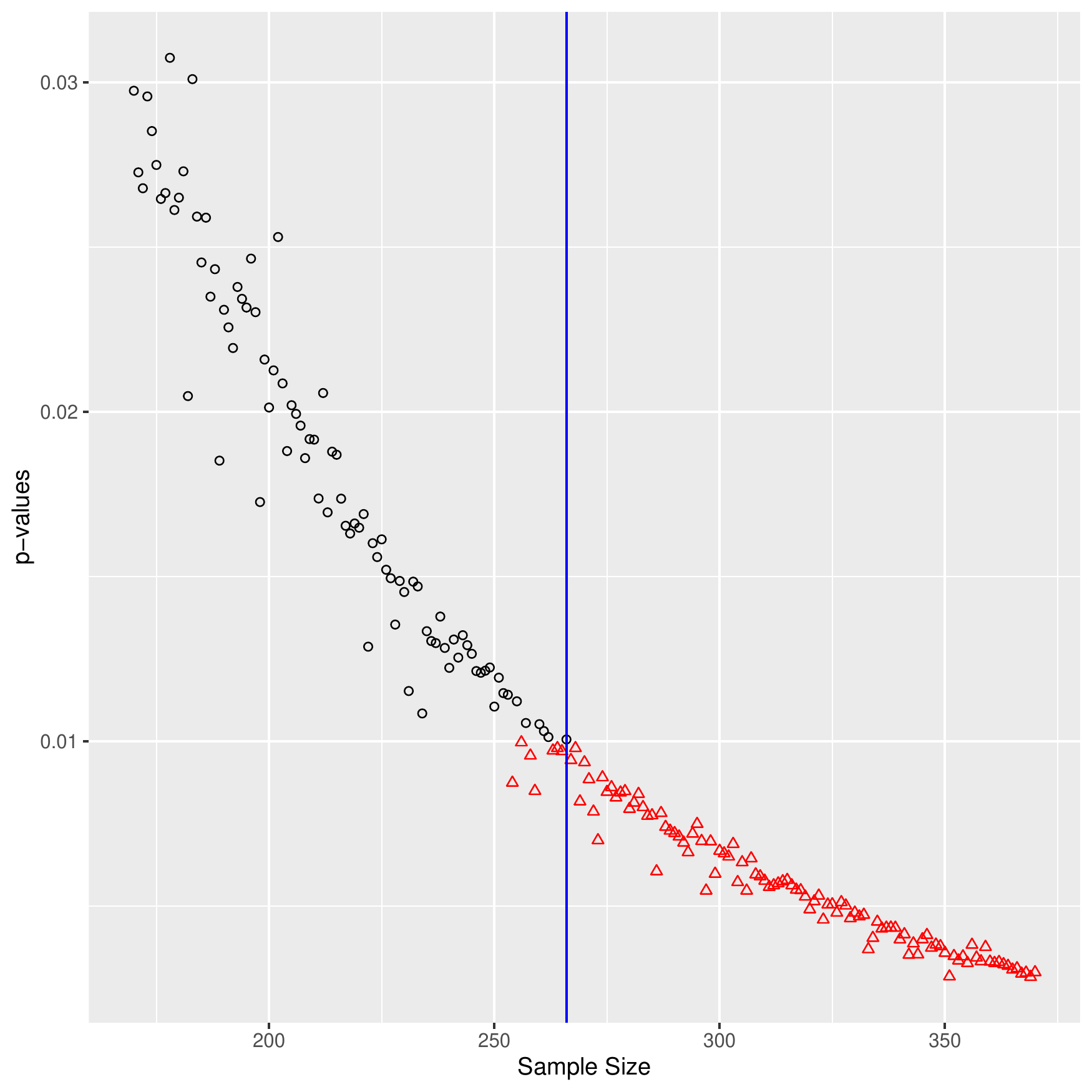}\\
			\caption{\small{p-values for different sample sizes in the 2-spin Curie-Weiss model ($\beta_0= 0.7>\beta^*(2)$ and $\beta= 0.9$).}}
			\label{scplot1}
		\end{minipage}%
		\begin{minipage}{0.5\textwidth}
			
			
			\centering
			\includegraphics[width=1\linewidth, height=0.28\textheight]
			{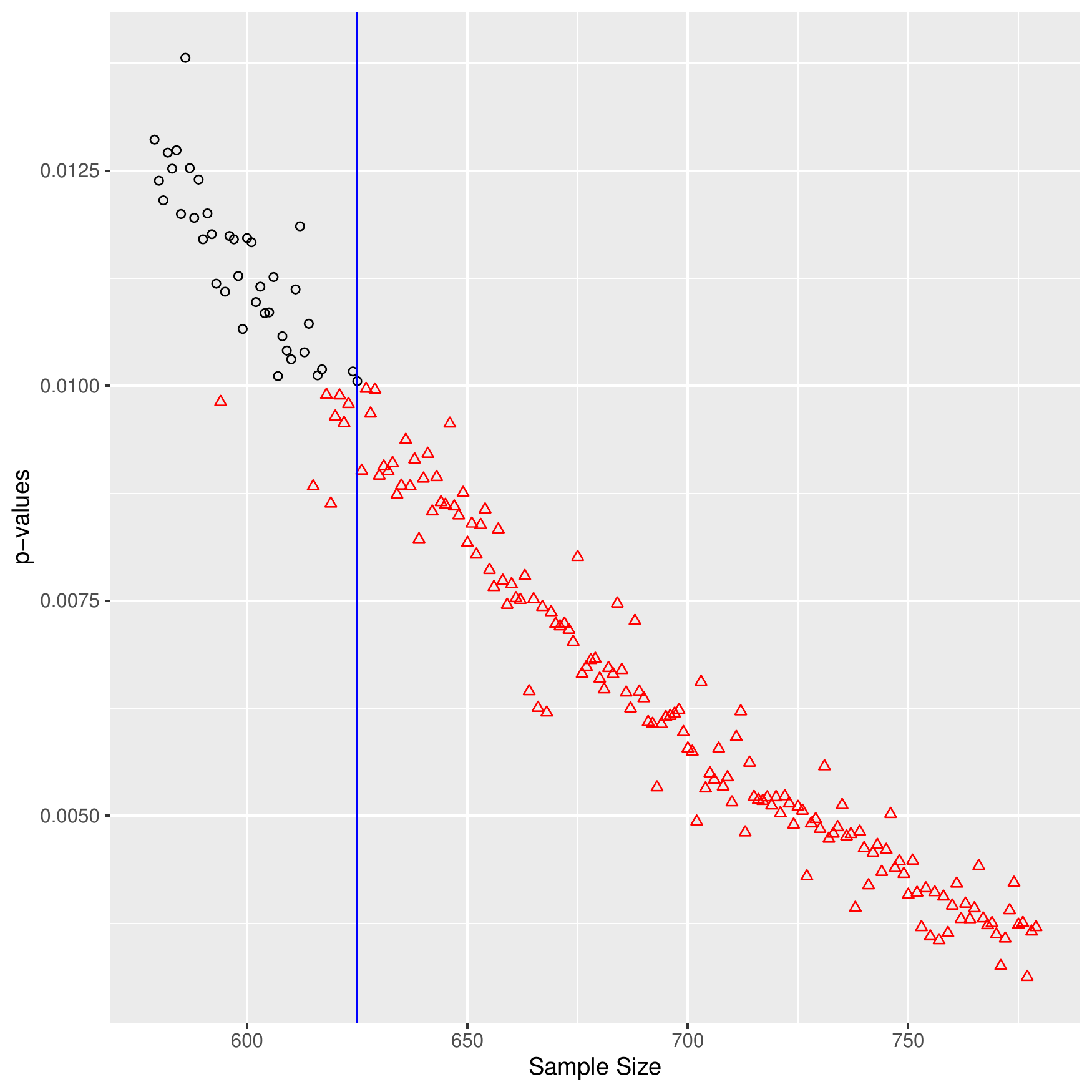}\\
			\caption{\small{p-values for different sample sizes in the 3-tensor Curie-Weiss model ($\beta_0= 0.68\in (\beta^*(3),\log 2),$ $\beta= 0.9$).}}
			\label{scplot2}
		\end{minipage}
	\end{figure}

	Now, we illustrate our theoretical results with numericals obtained from simulated data.  In Figure \ref{scplot1}, we plot the average $p$-value of the MPL test obtained from $10,000$ tuples generated from the $2$-spin Curie-Weiss model at $\beta = 0.9$ against the null $\beta_0 = 0.7 > \beta^*(2)$, for each sample size ranging from around $175$ to $375$. We see that from around sample size $266$, the average $p$-value goes down (and remains) below $\delta := 0.01$. This matches very closely with the theoretical sample-size value of $270$ that one will obtain in this setting, from Theorem \ref{main_mpl}. Figure \ref{scplot2} illustrates the average $p$-value of the MPL test obtained from $10,000$ tuples generated from the $3$-tensor Curie-Weiss model at $\beta = 0.9$ against the null $\beta_0 = 0.68 \in (\beta^*(3),\log 2)$, for each sample size ranging from around $575$ to $775$. We see that from around sample size $625$, the average $p$-value goes down (and remains) below $\delta := 0.01$. Once again, this matches somewhat closely with the theoretical asymptotic sample-size value of $679$ for the MPL test that one will obtain in this setting, from Theorem \ref{main_mpl}. In this case, the theoretical asymptotic sample size for the ML test turns out to be $533$. This is smaller than the theoretical and empirically obtained sample sizes of $679$ and $625$ (respectively) for the MPL test, thereby demonstrating our theoretical finding that the MPL test becomes much less efficient than the ML test for $\beta_0 < \log 2$ and sufficiently high $\beta$. 
	
\begin{figure*}[h]\vspace{0in}
		\centering
		\begin{minipage}[l]{1.0\textwidth}
			\centering
			\includegraphics[height=3.9in,width=5in]
			{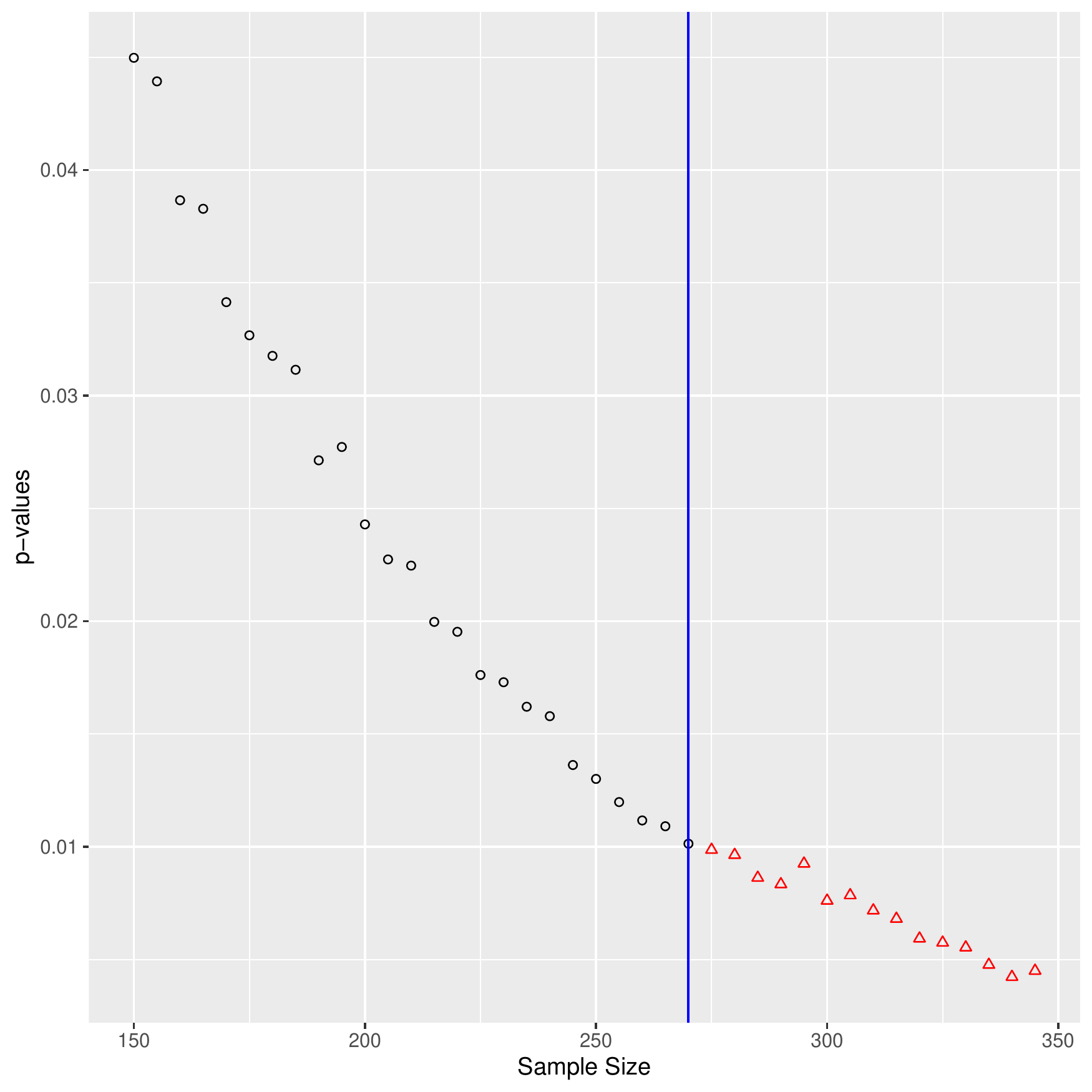}\\
		\end{minipage}
		\caption{\small{p-values for different sample sizes in the 2-spin \ern~model ($\beta_0= 0.7>\beta^*(2)$ and $\beta= 0.9$).}}
		\label{scplot3}
	\end{figure*}

	Next, we consider the Erd\H{o}s-R\'enyi Ising model in our numerical studies. The strategy is to generate a random matrix $\bm A$ with independent Bernoulli$(\alpha_n)$ entries, and simulate $1000$ samples under both the null and alternative Erd\H{o}s-R\'enyi Ising distributions using Glauber dynamics. We vary $n$ in steps of size $5$ in a region around the asymptotic optimal sample size given by Theorem~\ref{ermain}. Given $\boldsymbol{A}$, we generate each sample using $10^6$ iterations of the Glauber dynamics repeating the simulation independently 1000 times each under the null distribution $\beta_0$ and the alternative distribution $\beta \neq \beta_0$.
	 In Figure \ref{scplot3}, we plot the average $p$-values of the MPL tests obtained from $1000$ tuples generated from the $2$-spin \ern~model at $\beta_0 = 0.7 > \beta^*(2)$ against the alternative $\beta = 0.9$, for sample size ranging from around $150$ to $350$.
The simulation results match with the theoretical results indicating that the chains in the Glauber dynamics have mixed properly. 

	\section{Proofs of the Main Results}\label{mainproofr}
	In this section, we prove the main results stated in Sections \ref{sec2} and \ref{sec:ermodel}.  We start with the proof of Theorem \ref{main_mpl}.
	\subsection{Proof of Theorem \ref{main_mpl}}\label{mainmplproof7} 
	A main ingredient in the proof of Theorem \ref{main_mpl} is the verification of conditions \eqref{i1} and \eqref{i2} in \cite{bahadur} . We begin with the verification of condition \eqref{i1} with almost sure convergence replaced by convergence in probability.
	
	\begin{lemma}\label{cond1}
		Under every $\beta > \beta^*(p)$, we have:
		$$n^{-1/2} T_n\xrightarrow{P} \beta-\beta_0\quad\textrm{as}~n\rightarrow \infty~,$$
		where $T_n$ is either $\sqrt{n}(\bmp -\beta_0)$ or $\sqrt{n}(\bml -\beta_0)$.
	\end{lemma}

	We will now verify condition \ref{i2}. For this, we will need the following lemma on the large deviation of $\os$, which follows from \cite{pspinldp}. 
	
	\begin{lemma}\label{xbarldp}
		For every subset $A \subseteq [-1,1]$ such that $A^o$ is dense in $\overline{A}$, we have:
		$$\lim_{n\rightarrow \infty} \frac{1}{n} \log \p_{\beta,p}(\os \in A) = \sup_{x \in A} H_{\beta,p}(x) - \sup_{x\in [-1,1]} H_{\beta,p}(x)~.$$
		where $H_{\beta,p}$ is as defined in \eqref{hdef}.
	\end{lemma}
	
	\noindent Lemmas \ref{cond1} and \ref{xbarldp} are proved in Appendix \ref{maintr}. We are now ready to prove Theorem \ref{main_mpl}.
	\begin{proof}[Proof of Theorem \ref{main_mpl}]
		We deal with the test based on the MPL estimator first. To begin with,  note that $\bmp = \eta_p(\os)$. Fix $t>0$, whence we have by Lemma \ref{xbarldp}:
		\begin{eqnarray*}
			n^{-1}\log[1-F_{T_n,\beta_0}(n^{1/2} t)] &=& n^{-1} \log \p_{\beta_0,p}\left(n^{-1/2}T_n > t\right)\\&=& n^{-1} \log \p_{\beta_0,p}\left(\bmp-\beta_0>t\right)\\&=& n^{-1}\log \p_{\beta_0,p}\left(\eta_p(\os) > \beta_0 + t\right)\\&=& n^{-1} \log \p_{\beta_0,p}\left(\os \in \eta_p^{-1} ((\beta_0+t,\infty))\right)\\&=& \sup_{x \in \eta_p^{-1} ((\beta_0+t,\infty))} H_{\beta_0,p}(x) - \sup_{x\in [-1,1]} H_{\beta_0,p}(x) +o(1).
		\end{eqnarray*}
		where the last step follows from Lemma \ref{xbarldp}, since it follows from the proof of Lemma \ref{mlmplequal}, that the set $\eta_p^{-1} ((\beta_0+t,\infty))$ is a union of finitely many disjoint, non-degenerate intervals, and hence, its interior is dense in its closure. 
		
		In view of the above discussion, we conclude that the function $f$ in condition \eqref{i2} is given by:
		
		$$f(t) = \sup_{x\in [-1,1]} H_{\beta_0,p}(x) - \sup_{x \in \eta_p^{-1} ((\beta_0+t,\infty))} H_{\beta_0,p}(x)~.$$
		Since $x\mapsto H_{\beta_0,p}(x)$ and $\beta \mapsto m_*(\beta,p)$ are continuous functions (by Lemma \ref{mprop}) on $[-1,1]$ and $(\beta^*(p),\infty)$ respectively, we conclude (in view of Lemma \ref{mlmplequal}) that $f$ is continuous on an open neighborhood of $\beta-\beta_0$. Also, in view of Lemma \ref{mlmplequal}, the argument given below for the ML estimator, and the fact that $\beta_0>\beta^*(p)$, it will follow that $f>0$ on a non-empty open neighborhood of $\beta-\beta_0$.
		
		The Bahadur slope of $\bmp$ at an alternative $\beta$ is then given by $2f(\beta-\beta_0)$. This completes the proof of Theorem \ref{main_mpl} for the test based on the MPL estimator.
		
		For the test based on the ML estimator, note that for every $t > 0$, we have:
		\begin{eqnarray*}
			n^{-1}\log[1-F_{T_n,\beta_0}(n^{1/2} t)] &=& n^{-1} \log \p_{\beta_0,p}\left(n^{-1/2}T_n > t\right)\\&=& n^{-1} \log \p_{\beta_0,p}\left(\bml-\beta_0>t\right)\\&=& n^{-1}\log \p_{\beta_0,p}\left(\os^p > \e_{\beta_0 + t,p}(\os^p)\right)~.
		\end{eqnarray*}
		The last step follows from the following facts: 
		\begin{enumerate}
			\item The function $F_n(\beta,p) := \log Z_n(\beta,p)$ is strictly convex in $\beta$ (Lemma C.5 in \cite{jaesungcw}) and hence, $\frac{\partial F_N(\beta,p)}{\partial \beta}$ is strictly increasing in $\beta$;
			\item The ML equation is given by $\frac{\partial F_N(\beta,p)}{\partial \beta}\Big|_{\beta = \bml} = N \os^p ~;$
			\item $\frac{\partial F_N(\beta,p)}{\partial \beta} = \e_{\beta,p} (N \os^p)$.
		\end{enumerate}
		Now, it follows from \cite{jaesungcw} and the dominated convergence theorem, that $$\e_{\beta_0+t,p}(\os^p) \rightarrow m_*(\beta_0+t,p)^p~.$$ Fix $\varepsilon \in (0,m_*(\beta_0+t,p))$, to begin with. Then, there exists $N\ge 1$, such that $$\left(m_*(\beta_0+t,p)-\varepsilon\right)^p \le \e_{\beta_0+t,p}(\os^p) \le \left(m_*(\beta_0+t,p)+\varepsilon\right)^p$$ for all $n \ge N$. Let us first consider the case $p$ is odd. We then have the following by Lemma \ref{xbarldp}:
		\begin{eqnarray*}
			\limsup_{n\rightarrow \infty} n^{-1}\log[1-F_{T_n,\beta_0}(n^{1/2} t)] &=& \limsup_{n\rightarrow \infty} n^{-1}\log \p_{\beta_0,p}\left(\os^p > \e_{\beta_0 + t,p}(\os^p)\right)\\&\le& \limsup_{n\rightarrow \infty} n^{-1} \log \p_{\beta_0,p} \left(\os^p > \left(m_*(\beta_0+t,p)-\varepsilon\right)^p\right)\\&=& \limsup_{n\rightarrow \infty} n^{-1} \log \p_{\beta_0,p} \left(\os > m_*(\beta_0+t,p)-\varepsilon\right)\\ &=& \sup_{x~>~m_*(\beta_0+t,p)-\varepsilon} H_{\beta_0,p}(x) - \sup_{x\in [-1,1]} H_{\beta_0,p}(x)~. 
		\end{eqnarray*}
		Similarly, we also have:
		\begin{eqnarray*}
			\liminf_{n\rightarrow \infty} n^{-1}\log[1-F_{T_n,\beta_0}(n^{1/2} t)] &=& \liminf_{n\rightarrow \infty} n^{-1}\log \p_{\beta_0,p}\left(\os^p > \e_{\beta_0 + t,p}(\os^p)\right)\\&\ge& \liminf_{n\rightarrow \infty} n^{-1} \log \p_{\beta_0,p} \left(\os^p > \left(m_*(\beta_0+t,p)+\varepsilon\right)^p\right)\\&=& \liminf_{n\rightarrow \infty} n^{-1} \log \p_{\beta_0,p} \left(\os > m_*(\beta_0+t,p)+\varepsilon\right)\\ &=& \sup_{x~>~m_*(\beta_0+t,p)+\varepsilon} H_{\beta_0,p}(x) - \sup_{x\in [-1,1]} H_{\beta_0,p}(x)~. 
		\end{eqnarray*}
		Since $\varepsilon > 0$ can be arbitrarily small, and $H_{\beta,p}$ is continuous, we must have for all odd $p$:
		\begin{equation}\label{oddp1}
			\lim_{n\rightarrow \infty} n^{-1}\log[1-F_{T_n,\beta_0}(n^{1/2} t)] = \sup_{x ~> ~m_*(\beta_0+t,p)} H_{\beta_0,p}(x) - \sup_{x\in [-1,1]} H_{\beta_0,p}(x)~.
		\end{equation}
		
		Next, suppose that $p$ is even. In this case, $\bm X$ and $-\bm X$ have the same distribution, and hence, so do $\os$ and $- \os$. Hence, for every positive real number $\alpha$, we have:
		\begin{eqnarray*}
			n^{-1} \log \p_{\beta_0,p} \left(\os^p > \alpha^p\right) &=& n^{-1} \log \left[2~\p_{\beta_0,p} \left(\os > \alpha\right)\right] \\&=& n^{-1}\log \p_{\beta_0,p}(\os >\alpha) + o(1)~.
		\end{eqnarray*}
		Hence, the same argument as for the case of odd $p$ also works here, showing that \eqref{oddp1} holds when $p$ is even, too.
		
		In view of the above discussion, we conclude that the function $f$ in condition \eqref{i2} is given by:
		
		$$f(t) = \sup_{x\in [-1,1]} H_{\beta_0,p}(x) - \sup_{x ~>~ m_*(\beta_0+t,p)} H_{\beta_0,p}(x)~.$$
		Since $x\mapsto H_{\beta_0,p}(x)$ and $\beta \mapsto m_*(\beta,p)$ are continuous functions (by Lemma \ref{mprop}) on $[-1,1]$ and $(\beta^*(p),\infty)$ respectively, we conclude that $f$ is continuous on an open neighborhood of $\beta-\beta_0$. Also, 
		$f(\beta-\beta_0)>0$ (and hence, $f(t) >0$ on a non-empty open neighborhood of $\beta-\beta_0$), since $H_{\beta_0,p}$ is strictly decreasing on $m_*(\beta,p)$, and $m_*(\beta,p) > m_*(\beta_0,p)$ (by Lemma \ref{mprop}).
		
		The Bahadur slope of $\bml$ at an alternative $\beta$ is thus given by $2f(\beta-\beta_0)$. This completes the proof of Theorem \ref{main_mpl} for the test based on the ML estimator. The proof of Theorem \ref{main_mpl} is now complete.
	\end{proof}
	
	\subsection{Proof of Theorem \ref{p2case}}\label{p2cproof}
	The result for $p=2$ follows directly from Theorem \ref{main_mpl} and the following lemma, which is proved in Appendix \ref{maintr}.
	
	\begin{lemma}\label{mlmplequal}
		For every $\beta > \beta_0 > \beta^*(p)$, we have:
		\[   
		\sup_{x \in \eta_p^{-1} ((\beta,\infty))} H_{\beta_0,p}(x) = 
		\begin{cases}
			\sup_{x ~>~ m_*(\beta,p)} H_{\beta_0,p}(x) &\quad\text{if} ~p=2,\\
			\max\left\{\sup_{x ~>~ m_*(\beta,p)} H_{\beta_0,p}(x)~,~0\right\} &\quad\text{if} ~p\ge 3.\\
		\end{cases}
		\] 
	\end{lemma}
	
	\subsection{Proof of Theorem \ref{p3case}}\label{p3casepr7}
	The proof of Theorem \ref{p3case} depends crucially on the following important lemma:
	
	\begin{lemma}\label{nesu}
		For $p\ge 3$, $\beta>\beta_0>\beta^*(p)$ and $\delta \in (0,1)$, a necessary and sufficient condition for $c_{\bmp}(\beta_0,\beta,p) = c_{\bml}(\beta_0,\beta,p) \iff N_{\bmp}^*(\beta_0,\beta,\delta,p) = N_{\bml}^*(\beta_0,\beta,\delta,p)$ is $ H_{\beta_0,p}(m_*(\beta,p)) \ge 0$.
	\end{lemma}
	
	\begin{proof}[Proof of Lemma \ref{nesu}]
		For $p\ge 3$, it follows from Lemma \ref{hpr} that $H_{\beta,p}'\le 0$ on $(m_*(\beta,p),1)$, and hence, $H_{\beta_0,p}'\le 0$ on $(m_*(\beta,p),1)$. Consequently, 
		\begin{equation}\label{eq2}
			\sup_{x ~>~ m_*(\beta,p)} H_{\beta_0,p}(x) = \hbop(m_*(\beta,p))~.
		\end{equation}
		Lemma \ref{nesu} now follows from Lemma \ref{mlmplequal}.
	\end{proof}
	
	Returning to the proof of Theorem \ref{p3case}, note that it follows from \eqref{eq2} that $$\hbop(m_*(\beta,p)) \ge \hbop(1) = \beta_0-\log2~.$$
	This shows that the condition $\beta_0\ge \log 2$ is sufficient to ensure equality of the Bahadur slopes and the asymptotic optimal sample sizes. On the other hand, if $\beta_0<\log 2$, then $\lim_{x\rightarrow 1} H_{\beta_0,p}(x) < 0$. Since $\lim_{\beta\rightarrow \infty} m_*(\beta,p)=1$ (by Lemma \ref{mprop}), we must have for all $\beta>\beta_0$ large enough,
	$H_{\beta_0,p}(m_*(\beta,p)) < 0$, which shows, in view of Lemma \ref{nesu}, that the Bahadur slopes and asymptotic optimal sample sizes for the tests based on the MPL and ML estimators do not agree in this case.
	
	\subsection{Proof of Theorem \ref{beta1fixed}}  \label{betafixedpr7}
	By Lemma \ref{nesu}, $\beta_0> \beta^*(p)$ satisfies $\mathrm{eff}(\bml,\bmp,\beta_0,\beta) > 1$ if and only if $H_{\beta_0,p}(m_*(\beta,p)) < 0$, i.e. $\beta_0 m_*(\beta,p)^p < I(m_*(\beta,p))$.
	This establishes the first part of Theorem \ref{beta1fixed}. Next, note that 
	\begin{equation}\label{den1eq}
		\lim_{\beta_0 \rightarrow \beta^*(p)^+} \sup_{x\in [-1,1]} H_{\beta_0,p}(x) = 0~.
	\end{equation}
	On the other hand, we also have:
	$$\lim_{\beta_0\rightarrow \beta^*(p)^+}\sup_{x~>~m_*(\beta,p)} H_{\beta_0,p}(x) = H_{\beta^*(p),p}(m_*(\beta,p)) < 0~.$$
	Hence, it follows from Lemma \ref{mlmplequal} that:
	\begin{equation}\label{den2eq}
		\lim_{\beta_0\rightarrow \beta^*(p)^+} \sup_{x \in \eta_p^{-1} ((\beta,\infty))} H_{\beta_0,p}(x) = \lim_{\beta_0\rightarrow \beta^*(p)^+}\max\left\{0~,~\sup_{x~>~m_*(\beta,p)} H_{\beta_0,p}(x)\right\} = 0~.
	\end{equation}
	It thus follows from \eqref{den1eq}, \eqref{den2eq} and Theorem \ref{main_mpl}, that 
	\begin{equation}\label{sasizeinf}
		\lim_{\beta_0\rightarrow \beta^*(p)^+} N_{\bmp}^*(\beta_0,\beta,\delta,p) = \infty \quad\textrm{and}\quad\lim_{\beta_0\rightarrow \beta^*(p)^+} N_{\bml}^*(\beta_0,\beta,\delta,p) = \frac{\log(\delta)}{H_{\beta^*(p),p}(m_*(\beta,p))} < \infty.
	\end{equation}   	
	Finally, the limiting Bahadur efficiency part follows directly from \eqref{sasizeinf}. This completes the proof of Theorem \ref{beta1fixed}.
	
	\subsection{Proof of Theorem \ref{stablesamp}}  \label{th5finpr}
	Since $m_*(\beta,p) \rightarrow 1$ as $\beta \rightarrow \infty$ (by Lemma \ref{mprop}), we have:
	$$\lim_{\beta\rightarrow \infty} H_{\beta_0,p}(m_*(\beta,p)) = \beta_0 - \log 2 < 0~.$$
	Hence, $H_{\beta_0,p}(m_*(\beta,p)) < 0$ for all $\beta$ large enough, which shows, in view of Lemma \ref{mlmplequal}, that:
	$$\sup_{x \in \eta_p^{-1} ((\beta,\infty))} H_{\beta_0,p}(x) = 0~.$$ Theorem \ref{stablesamp} now follows from Theorem \ref{main_mpl}.
	
	\subsection{Proof of Theorem \ref{ermain}}\label{erproof7}
	As in the Curie-Weiss model, a main ingredient in the proof of Theorem \ref{ermain} is the verification of conditions \eqref{i1} and \eqref{i2} in \cite{bahadur}. Once again, we begin with the verification of condition \eqref{i1} with almost sure convergence replaced by convergence in probability. The following result is an analogous version of Lemma \ref{cond1} for the model \eqref{ermodel}.
	
	\begin{lemma}\label{cond12}
		Under the model \eqref{ermodel}, for every $\beta > \beta^*(p)$, we have:
		$$n^{-1/2} T_n\xrightarrow{P} \beta-\beta_0\quad\textrm{as}~n\rightarrow \infty~,$$
		where $T_n$ is either $\sqrt{n}(\bmp -\beta_0)$ or $\sqrt{n}(\bml -\beta_0)$.
	\end{lemma}
	
	We will now verify condition \ref{i2}. Towards this, we will use Lemma \ref{hamiltcomp} to compare the probability models \eqref{cwmodel} and \eqref{ermodel}. In fact, we prove a slightly more general result below, which solves our purpose, but may be of independent interest for more general objectives.
	
	\begin{lemma}\label{measureclose}
		Let $\{\beta_n\}_{n \geq 1}$ be a bounded sequence of positive real numbers.
		Then, with probability $1$, we have: $$\sup_{A, B \subseteq \sa}\big|\log \p_{\beta_n,p}^*(A|B) - \log \p_{\beta_n,p}(A|B)\big| = O(n\gamma_n)~.$$
	\end{lemma} 
	
	\begin{remark}
		One can compare the logarithms of the two (unconditional) measures $\p^*$ \eqref{ermodel} and $\p$ \eqref{cwmodel} by taking $B := \{-1,1\}^n$ in Lemma \ref{measureclose}.
	\end{remark}
	
	One can now use Lemmas \ref{hamiltcomp} and \ref{measureclose} to compare the log-normalizing constants and asymptotics of the sample mean in the two models \eqref{cwmodel} and \eqref{ermodel}. 
	
	\begin{lemma}\label{freeenergy}
		We have the following with probability $1$.
		\begin{itemize}
			\item[1.]~If $Z_n^*(\beta,p)$ denotes the normalizing constant of the model \eqref{ermodel}, then $$|\log Z_n(\beta,p) - \log Z_n^*(\beta,p)| = O(n\gamma_n)~.$$
			\item[2.]~If $\bs$ is generated from the model \eqref{ermodel}, then for every $\beta > \beta^*(p)$ and fixed $\varepsilon > 0$, we have:
			\begin{equation}\label{wl}
				\p_{\beta,p}^*\left(|\os^p - m_*(\beta,p)^p| \ge \varepsilon\right) \le e^{-n\Omega(1)}~.
			\end{equation}
			In particular, $\os \xrightarrow{P} m_*(\beta,p)$ under the model \eqref{ermodel}.
		\end{itemize}
	\end{lemma}
	
	Lemmas \ref{measureclose} and \ref{freeenergy} are proved in Appendix \ref{maintr}. In view of all that we have above, we are now ready to prove Theorem \ref{ermain}.
	
	\noindent\textit{Proof of Theorem \ref{ermain}.} 
	We begin with the ML estimator first. Note that by Lemma \ref{hamiltcomp}, we have:
	\begin{eqnarray*}
		\p_{\beta_0,p}^*(\os^p > \e_{\beta_0+t,p}^*\os^p +6\gamma_n) &\le& \p_{\beta_0,p}^*(H_n(\bm X) > \e_{\beta_0+t,p}^*H_n(\bm X))\\ &\le& \p_{\beta_0,p}^*(\os^p > \e_{\beta_0+t,p}^*\os^p -6\gamma_n). 	
	\end{eqnarray*}
	Now, note that by part (2) of Lemma \ref{freeenergy} and the dominated convergence theorem, $\e_{\beta_0+t,p}^*\os^p \rightarrow m_*(\beta_0+t,p)^p$. Hence, we have:
	\begin{eqnarray*}
		\p_{\beta_0,p}^*(\os^p > m_*(\beta_0+t,p)^p +o(1)) &\le& \p_{\beta_0,p}^*(H_n(\bm X) > \e_{\beta_0+t,p}^*H_n(\bm X))\\ &\le& \p_{\beta_0,p}^*(\os^p > m_*(\beta_0+t,p)^p +\tilde{o}(1))
	\end{eqnarray*}
	where $o(1)$ and $\tilde{o}(1)$ denote two real sequences converging to $0$. Hence, for every fixed $\varepsilon>0$ sufficiently small, one has the following for all large $n$:
	\begin{eqnarray*}
		\p_{\beta_0,p}^*(\os^p > (m_*(\beta_0+t,p)+\varepsilon)^p) &\le& \p_{\beta_0,p}^*(H_n(\bm X) > \e_{\beta_0+t,p}^*H_n(\bm X))\\ &\le& \p_{\beta_0,p}^*(\os^p > (m_*(\beta_0+t,p)-\varepsilon)^p).
	\end{eqnarray*}
	It now follows from Lemma \ref{measureclose}, that:
	\begin{eqnarray*}
		\frac{1}{n}\log \p_{\beta_0,p}(\os^p > (m_*(\beta_0+t,p)+\varepsilon)^p) + o(1)&\le& \frac{1}{n} \log  \p_{\beta_0,p}^*(H_n(\bm X) > \e_{\beta_0+t,p}^*H_n(\bm X))\\ &\le& \frac{1}{n} \log \p_{\beta_0,p}(\os^p > (m_*(\beta_0+t,p)-\varepsilon)^p) + o(1).
	\end{eqnarray*}
	Theorem \ref{ermain} for the ML estimator now follows from the proof of Theorem \ref{main_mpl} and Lemma \ref{cond12}.
	
	Next, we consider the MPL estimator. It follows from (2.3) in \cite{jaesung1}, that $\bmp$ is the least solution of the equation (in $\beta$):
	$$H_n(\bm X) = \sum_{i=1}^n m_i(\bm X) \tanh(p\beta m_i(\bm X))~.$$
	Define $\psi_n(\beta):= n^{-1}\sum_{i=1}^n m_i(\bm X) \tanh(p\beta m_i(\bm X))$. Since with probability $1$, we have the following for all large $n$ $$\psi_n'(\beta) = \frac{p}{n}\sum_{i=1}^n m_i^2(\bm X) \mathrm{sech}^2(p\beta m_i(\bm X)) > 0,$$ the function $\psi_n$ is strictly increasing for all large $n$, with probability $1$. Hence, we have (with probability $1$ for all large $n$) the following for all $t>0$:
	
	\begin{eqnarray*}
		\frac{1}{n}\log \p_{\beta_0,p}^*(\bmp > \beta_0+t) &=& \frac{1}{n}\log \p_{\beta_0,p}^* (\psi_n(\bmp) > \psi_n(\beta_0+t))\\&=& \frac{1}{n}\log \p_{\beta_0,p}^* \left(\frac{1}{n} H_n(\bm X) > \frac{1}{n}\sum_{i=1}^n m_i(\bm X)\tanh(p(\beta_0+t) m_i(\bm X))\right)\\&=& \frac{1}{n}\log \p_{\beta_0,p}^* \left(\os^p > \frac{1}{n}\sum_{i=1}^n m_i(\bm X)\tanh(p(\beta_0+t) m_i(\bm X))+o(1)\right)
	\end{eqnarray*}
	
	Therefore, in view of Lemma \ref{hamiltcomp}, Corollary \ref{hcm1} and Lemma \ref{measureclose}, we have the following for every fixed $\varepsilon > 0$ sufficiently small:
	\begin{eqnarray*}
		&& \limsup_{n\rightarrow \infty} \frac{1}{n}\log \p_{\beta_0,p}^*(\bmp > \beta_0+t)\\&\le & \limsup_{n\rightarrow \infty} \frac{1}{n}\log \p_{\beta_0,p}^* \left(\os^p > \frac{1}{n}\sum_{i=1}^n \os^{p-1}\tanh(p(\beta_0+t) \os^{p-1}) +o(1)\right)\\&= & \limsup_{n\rightarrow \infty} \frac{1}{n}\log \p_{\beta_0,p}^* \left(\os^p > \os^{p-1}\tanh(p(\beta_0+t) \os^{p-1}) +o(1)\right)\\&\le& \limsup_{n\rightarrow \infty} \frac{1}{n}\log \p_{\beta_0,p}^* \left(\eta_p(\os) > \beta_0+t-\varepsilon\right)\\&=& \limsup_{n\rightarrow \infty} \frac{1}{n}\log \p_{\beta_0,p} \left(\eta_p(\os) > \beta_0+t-\varepsilon\right)\\&=& \sup_{x \in \eta_p^{-1} ((\beta_0+t-\varepsilon,\infty))} H_{\beta_0,p}(x) - \sup_{x\in [-1,1]} H_{\beta_0,p}(x).
	\end{eqnarray*}
	We can now take $\varepsilon\downarrow 0$ to conclude that:
	\begin{equation}\label{mpler1}
		\limsup_{n\rightarrow \infty} \frac{1}{n}\log \p_{\beta_0,p}^*(\bmp > \beta_0+t) \le \sup_{x \in \eta_p^{-1} ((\beta_0+t,\infty))} H_{\beta_0,p}(x) - \sup_{x\in [-1,1]} H_{\beta_0,p}(x)~.
	\end{equation}
	By an exactly similar approach, we can show that:
	\begin{equation}\label{mpler2}
		\liminf_{n\rightarrow \infty} \frac{1}{n}\log \p_{\beta_0,p}^*(\bmp > \beta_0+t) \ge \sup_{x \in \eta_p^{-1} ((\beta_0+t,\infty))} H_{\beta_0,p}(x) - \sup_{x\in [-1,1]} H_{\beta_0,p}(x)~.
	\end{equation}
	Theorem \ref{ermain} now follows from \eqref{mpler1} and \eqref{mpler2}. \qed

	\section{Discussion}\label{discussion}
	In this paper, we derived the Bahadur slopes and optimal sample sizes required for significance of the tests based on the maximum likelihood (ML) and the maximum pseudolikelihood (MPL) estimators for the tensor Curie-Weiss model. One of our interesting findings is that although the MPL estimator is just an approximation of the ML estimator, the former is as Bahadur efficient as the latter everywhere in the parameter space for $p=2$, and throughout most of the parameter space for $p\ge 3$. More precisely, the MPL estimator is equally Bahadur efficient as the ML estimator for $p=2$, and for $p\ge 3$ this is true for all values of the alternative, if and only if the null parameter is greater than or equal to $\log 2$. For $p\ge 3$, if the null parameter lies strictly between the estimation threshold and $\log 2$, then the MPL estimator is strictly less efficient than the ML estimator for all sufficiently large values of the alternative parameter. However, even in this regime, the Bahadur ARE of the MPL estimator with respect to the ML estimator is lower bounded by a positive fraction. We also showed that our results hold verbatim for the more general class of tensor Erd\H{o}s-R\'enyi Ising models, where we can even allow for some sparsity in the underlying hypergraph. 
	
	We conjecture that similar results are also true for Ising models on dense deterministic and stochastic block model hypergraphs (with the Hamiltonian being suitably scaled) too, and believe that this can be shown by slight (routine) modifications of the methods used in Section \ref{sec:ermodel}. A potentially interesting direction for future research in this area, would be to consider Ising models on more general hypergraphs, for example arbitrary regular hypergraphs. Probability limits and fluctuations of the average magnetization for $2$-spin Ising models on $d$-regular graphs have been recently derived in \cite{nabarun}, where the authors show that the fluctuations are universal and same as that of the $2$-spin Curie-Weiss model in the entire ferromagnetic parameter regime as long as $d \gg \sqrt{n}$. The next natural step would thus, be to derive a large deviation principle for the Hamiltonians of such models, not just in the $2$-spin case, but also for the tensor case. These will in turn enable one to derive the Bahadur slopes and the optimal sample sizes for the ML and the MPL estimators in Ising models on regular hypergraphs.
	
	Since the behavior of the MPL estimator in terms of its Bahadur ARE with respect to the ML estimator is dependent on whether $p=2$ or $p>2$, one may be interested in asking if it is possible to test the hypothesis $p=2$ versus $p>2$, at least in the Curie-Weiss model, based on a single observation $\bm X$ from the model. This is an interesting problem in its own right, and one plausible approach may be to estimate $\beta$ assuming $p=2$ (by either the ML or the MPL approach), simulate a large number $B$ of observations $\bm X^{(1)},\bm X^{(2)},\ldots,\bm X^{(B)}$ from the $2$-spin Curie-Weiss model with parameter $\hat{\beta}$, and determine whether $\overline{X}$ lies within the $2.5\%$ and $97.5\%$ quantiles of the histogram formed by $\overline{X^{(1)}},\overline{X^{(2)}},\ldots,\overline{X^{(B)}}$. If not, then that should be a reasonable evidence against the null hypothesis $p=2$. This approach has some flaws though, one being that we are never sure if the actual observation $\bm X$ is coming from a Curie-Weiss model at low temperature ($\beta$ above the threshold), because otherwise, $\hat{\beta}$ will be inconsistent. The second drawback is that one can never say surely that the observation is coming from a $2$-spin Curie-Weiss model indeed, even if $\overline{X}$ does lie within the $2.5\%$ and $97.5\%$ quantiles of the empirical histogram. In other words, the power guarantee of this testing approach is dubious. We leave this problem as an interesting direction for future research.  
	
	Finally, we would like to mention that although this paper is concerned with efficient estimation and testing of the coupling strength in tensor Ising models and the model considered here does not have any external magnetic field term, there is a significant literature on testing for external fields as well (see for example, \cite{rm1,rm2}), at least in the classical $2$-spin setting. In the Curie-Weiss model with external field, the sufficient statistic is still the sample mean $\overline{X}_n$, and both the ML and the MPL estimates are once again functions of $\overline{X}_n$. Hence, we expect similar techniques involving asymptotics and large deviations of $\overline{X}_n$ to apply in deriving the Bahadur slopes and optimal sample sizes corresponding to the estimates of the external fields, as well.
	
	\section{Acknowledgment}
	\noindent Somabha Mukherjee was supported by the National University of Singapore Start-Up Grant R-155-000-233-133, 2021. Swarnadip Ghosh was supported by the National Science Foundation BIGDATA grant IIS-1837931. The authors thank Bhaswar B. Bhattacharya and Paul Switzer for several helpful discussions and careful comments.

	\appendix
	
	\section{Proofs of Technical Lemmas}\label{maintr}
	In this section, we prove the technical lemmas mentioned in Sections \ref{sec2} and \ref{sec:ermodel}. 
	\subsection{Proof of Lemma \ref{cond1}}
	Note that $n^{-1/2}T_n = \hat{\beta}-\beta_0$, where $\hat{\beta}$ is either $\bml$ or $\bmp$. It follows from \cite{jaesungcw} and \cite{jaesung1}, that under every $\beta > \beta^*(p)$, $\hat{\beta}\xrightarrow{P} \beta$. This proves Lemma \ref{cond1}.
	
	\subsection{Proof of Lemma \ref{xbarldp}}
	It follows from display (18) in \cite{pspinldp} that $\os$ satisfies a large deviation principle (LDP) with rate function:
	$$I(x) := - \beta x^p + \frac{x}{2} \sinh^{-1}\left(\frac{2x}{1-x^2}\right) + \frac{1}{2} \log\left(1-x^2\right) - \inf_{y\in \mathbb{R}} \left\{\sup_{z\in \mathbb{R}} \{yz - \log \cosh(z)\} - \beta y^p  \right\}~.$$
	Using the identity $$\sinh^{-1}(x) = \log\left(x + \sqrt{x^2+1}\right),$$
	we have:
	\begin{eqnarray*}
		&&  - \beta x^p + \frac{x}{2} \sinh^{-1}\left(\frac{2x}{1-x^2}\right) + \frac{1}{2} \log\left(1-x^2\right)\\&=& - \beta x^p + \frac{x}{2} \log\left(\frac{2x}{1-x^2} + \sqrt{\frac{4x^2}{(1-x^2)^2} +1}\right) + \frac{1}{2} \log\left(1-x^2\right)\\&=& - \beta x^p + \frac{x}{2} \log\left(\frac{(1+x)^2}{1-x^2} \right) + \frac{1}{2} \log\left(1-x^2\right)\\&=& -\beta x^p + x \log(1+x) + \frac{1-x}{2} \log(1-x^2)\\&=& -\beta x^p + x \log(1+x) + \frac{1-x}{2} \log(1+x) + \frac{1-x}{2}\log(1-x)\\&=& - \beta x^p + \frac{1}{2}\left\{(1+x)\log(1+x) + (1-x)\log(1-x)\right\} = -H_{\beta,p}(x)~.
	\end{eqnarray*}
	It now follows from Lemma \ref{sup} that
	\begin{equation}\label{iexp}
		I(x) = -H_{\beta,p}(x) + \sup_{y \in [-1,1]} H_{\beta,p}(y)~.
	\end{equation}
	Lemma \ref{xbarldp} now follows from \eqref{iexp} and the fact that $I$ is a continuous function on $[-1,1]$.
	
	\subsection{Proof of Lemma \ref{mlmplequal}}
	First, suppose that $p$ is even. Then, $\eta_p$ is an even function, and hence, the set $\eta_p^{-1}((\beta,\infty))$ is symmetric around $0$. This, together with the fact that $H_{\beta_0,p}$ is an even function, implies that
	\begin{equation}\label{rc}
		\sup_{x \in \eta_p^{-1} ((\beta,\infty))} H_{\beta_0,p}(x) = \sup_{x \in \eta_p^{-1} ((\beta,\infty))\bigcap(0,1]} H_{\beta_0,p}(x)~.
	\end{equation}  
	Now, note that if $p$ is odd, then $\eta_p^{-1}((\beta,\infty)) \bigcap [-1,0] =\emptyset$, since $\eta_p(x) \le 0$ for all $x \in [-1,0]$. Hence, \eqref{rc} is valid for odd $p$, too.
	
	Now, $x \in \eta_p^{-1} ((\beta,\infty))~\bigcap ~(0,1]$ if and only if $x \in (0,1]$ satisfies $p^{-1}x^{1-p}\tanh^{-1}(x) > \beta$, if and only if $x \in (0,1]$ satisfies $$H_{\beta,p}'(x) =\beta p x^{p-1} -\tanh^{-1}(x) < 0~.$$
	If $p\ge 3$, then by Lemma \ref{hpr}, this region is precisely equal to $\left\{(0,\um(\beta,p))~\bigcup ~(m_*(\beta,p), 1]\right\} \setminus F$ for some finite set $F$ (which is either singleton or empty). Hence, for $p \ge 3$, we have by continuity of $\hbop$, that:
	
	\begin{equation}\label{p1}
		\sup_{x \in \eta_p^{-1} ((\beta,\infty))} H_{\beta_0,p}(x) = \max\left\{\sup_{x\in (0,\um(\beta,p))} H_{\beta_0,p}(x) ~,~\sup_{x\in (m_*(\beta,p),1]} H_{\beta_0,p}(x)\right\}.
	\end{equation}
	Since $\hbp(0) = 0$ and $\hbp$ is decreasing on $(0,\um(\beta,p))$, we must have:
	$$\sup_{x\in (0,\um(\beta,p))} H_{\beta,p}(x)=0~.$$ Further, since $\hbop(0) = 0$ and $H_{\beta_0,p} \le H_{\beta,p}$ on $[0,1]$, we must have:
	$$\sup_{x\in (0,\um(\beta,p))} H_{\beta_0,p}(x)=0~.$$
	Lemma \ref{mlmplequal} for $p\ge 3$ now follows from \eqref{p1}. Now, let $p=2$. Then, $$\eta_p^{-1} ((\beta,\infty))~\bigcap ~(0,1] = (m_*(\beta,p),1]~.$$ Hence,
	$$\sup_{x \in \eta_p^{-1} ((\beta,\infty))} H_{\beta_0,p}(x) = \sup_{x ~>~ m_*(\beta,p)} H_{\beta_0,p}(x)~.$$
	This completes the proof of Lemma \ref{mlmplequal}.  
	
	\subsection{Proof of Lemma \ref{cond12}}
	Note that $n^{-1/2}T_n = \hat{\beta}-\beta_0$, where $\hat{\beta}$ is either $\bml$ or $\bmp$. All the following arguments are on the following event, which has probability $1$ in view of Lemma \ref{hamiltcomp}:
	$$\mathcal{E}:= \left\{\sup_{\bt \in \sa} \left|H_n(\bt) - n \oss^p\right| \le 3n\gamma_n~\text{for all but finitely many}~n\right\}~.$$
	Let us first consider the case $\hat{\beta}=\bml$. Then, for every fixed $t>0$, we have:
	$$\p_{\beta,p}^*(\bml > \beta + t) = \p_{\beta,p}^*(H_n(\bm X) > \e_{\beta+t,p}^* H_n(\bm X)) \le \p_{\beta,p}^*\left(\os^p > \e_{\beta+t,p}^* (\os^p) - 6\gamma_n\right)~.$$ Now, by part (2) of Lemma \ref{freeenergy} and the dominated convergence theorem, we have:
	$$\e_{\beta+t,p}^*(\os^p) \rightarrow m_*(\beta+t,p)^p~.$$ Hence, once again by part (2) of Lemma \ref{freeenergy}, we have:
	$$\p_{\beta,p}^*(\bml > \beta + t) \le \p_{\beta,p}^*(\os^p >  m_*(\beta+t,p)^p - o(1)) \le \p_{\beta,p}^* \left(|\os^p-m_*(\beta,p)^p| > \Omega(1)\right) =o(1).$$ 
	Similarly, we can show that for every $t \in (0,\beta-\beta^*(p))$,
	$$\p_{\beta,p}^*(\bml < \beta - t) \le \p_{\beta,p}^*(\os^p <  m_*(\beta-t,p)^p + o(1)) \le \p_{\beta,p}^* \left(|\os^p-m_*(\beta,p)^p| > \Omega(1)\right) =o(1).$$ Hence, we conclude that $\bml \xrightarrow{P} \beta$ under the model $\p_{\beta,p}^*$. This proves Lemma \ref{cond12} when $\hat{\beta}= \bml$.
	
	Now, suppose that $\hat{\beta} = \bmp$. By part 1 of Lemma \ref{freeenergy}, we have:
	$$\frac{1}{n} \log Z_n^*(\beta,p) = \frac{1}{n} \log Z_n(\beta,p) + o(1)~.$$
	By Theorem 2.3 in \cite{jaesung1} and Lemma \ref{c11}, we conclude that $\bmp$ is a consistent estimator of $\beta$ under the model $\p_{\beta,p}^*$. This completes the proof of Lemma \ref{cond12}. 
	
	\subsection{Proof of Lemma \ref{measureclose}}
	For any two sets $A, B \subseteq \sa$, note that:
	\begin{eqnarray}\label{conex}
		&&\p_{\beta_n,p}^*(A|B)\nonumber\\ &=& \frac{\sum_{\bt \in A\bigcap B} \exp\{\beta_n H_n(\bt)\}}{\sum_{\bt \in B} \exp\{\beta_n H_n(\bt) \}}\nonumber\\&=& \frac{\sum_{\bt \in A\bigcap B} \exp\{\beta_n \e H_n(\bt)\} e^{\beta_n (H_n(\bm x) - \e H_n(\bm x))}}{\sum_{\bt \in B} \exp\{\beta_n \e H_n(\bt)\} e^{\beta_n (H_n(\bm x) - \e H_n(\bm x)) } }~.
	\end{eqnarray}
	It follows from \eqref{conex} and Lemma \ref{hamiltcomp}, that with probability $1$, we have the following for all large $n$,
	\begin{equation}\label{conex2}
		e^{-6n\beta_n\gamma_n} \p_{\beta_n,p}(A|B)\leq \p_{\beta_n,p}^*(A|B) \leq e^{6n\beta_n\gamma_n} \p_{\beta_n,p}(A|B)~.
	\end{equation}
	Lemma \ref{measureclose} now follows on taking logarithm on both sides of \eqref{conex}, and recalling that the sequence $\{\beta_n\}_{n\geq 1}$ is bounded.
	
	\subsection{Proof of Lemma \ref{freeenergy}}
	\hfill\break
	\noindent\textit{Proof of} 1.~ Using Lemma \ref{hamiltcomp}, we have the following for all large $n$, with probability $1$:
	\begin{equation}\label{comppart}
		e^{-3\beta n\gamma_n}Z_n(\beta,p) \leq Z_n^*(\beta,p) \leq e^{3\beta n \gamma_n} Z_n(\beta,p)~.
	\end{equation}
	Part 1 follows on taking logarithm on both sides of \eqref{comppart}.
	\vskip.1in
	\noindent\textit{Proof of} 2.~ To begin with, define:
	\[   
	M_\varepsilon := 
	\begin{cases}
		(m_*(\beta,p)-\varepsilon, m_*(\beta,p) + \varepsilon) &\quad\text{if}~p~\text{is odd}\\
		(m_*(\beta,p)-\varepsilon, m_*(\beta,p) + \varepsilon) \bigcup (-m_*(\beta,p)-\varepsilon, -m_*(\beta,p) + \varepsilon) &\quad\text{if}~p~\text{is even}\\ 
	\end{cases}
	\]
	It follows from the arguments used in the proof of Lemma 3.1 in \cite{jaesungcw}, that
	\begin{equation}\label{mlemod}
		\p_{\beta,p}(\os \in M_\varepsilon^c) = O(n^\frac{3}{2}) \exp\left\{n\left(\sup_{t\in M_\varepsilon^c} H_{\beta,p}(t) - H_{\beta,p}(m_*(\beta,p))\right)   \right\} = e^{-n\Omega(1)} .
	\end{equation}
	It follows from \eqref{mlemod} and Lemma \ref{measureclose}, that 
	\begin{equation}\label{mlemod2}
		\p_{\beta,p}^*(\os \in M_\varepsilon^c) \le \p_{\beta,p}(\os \in M_\varepsilon^c) e^{O(n\gamma_n)} = e^{-n\Omega(1)} .
	\end{equation}
	Hence, we have:
	$$\p_{\beta,p}^*\left(|\os^p-m_*(\beta,p)^p| \ge \varepsilon\right) \le \p_{\beta,p}^*(\os \in M_{\varepsilon/p}^c) = e^{-n\Omega(1)}~.$$
	This completes the proof of Lemma \ref{freeenergy}.

	\section{Technical Results for the Tensor Curie-Weiss Model}\label{tech}
	In this section, we prove some technical results, which will be useful in the proofs of the main results concerning the tensor Curie-Weiss model.
	
	%
	%
	
	\begin{lemma}\label{sup}
		We have the following:
		\[   
		\sup_{z\in \mathbb{R}} \{yz - \log \cosh(z)\} = 
		\begin{cases}
			\frac{1}{2}\left\{(1+y)\log(1+y) + (1-y)\log(1-y)\right\} &\quad\text{if} ~y \in [-1,1],\\
			\infty &\quad\text{otherwise}~.\\
		\end{cases}
		\]
	\end{lemma}
	\begin{proof}
		Fix $y \in \mathbb{R}$ and define $g(z) := yz-\log \cosh(z)$. Let us begin with the case $y \in (-1,1)$. In this case, $g''(z) = -\mathrm{sech}^2(z) < 0$ and hence, $g$ is a strictly concave function. Consequently, any stationary point of $g$ is the unique global maximum of $g$. Since $g'(z) = y - \tanh(z)$, it follows that the only stationary point of $g$ is $\tanh^{-1}(y)$, and hence, 
		$$\sup_{z\in \mathbb{R}} \{yz - \log \cosh(z)\} = y \tanh^{-1}y - \log \cosh (\tanh^{-1}(y)) = y \tanh^{-1}(y) + \frac{1}{2} \log(1-y^2)~,$$
		where in the last step, we used the identity:
		$$\cosh (\tanh^{-1}(y)) = \frac{1}{\sqrt{1-y^2}}\quad\quad\textrm{for}~y\in (-1,1).$$
		The proof for the case $y \in (-1,1)$ now follows from the observation that
		$$y \tanh^{-1}(y) + \frac{1}{2} \log(1-y^2) = \frac{1}{2}\left\{(1+y)\log(1+y) + (1-y)\log(1-y)\right\}.$$
		
		Now, suppose that $y \ge 1$. Then, $g'(z) > 0$ for all $z \in \mathbb{R}$, and hence,
		$$\sup_{z\in \mathbb{R}} g(z) = \lim_{z\rightarrow \infty} g(z)~.$$
		Now, note that
		\[   
		\lim_{z\rightarrow \infty} e^{g(z)} = \lim_{z\rightarrow \infty} \frac{2e^{yz}}{e^z + e^{-z}} = \lim_{z\rightarrow \infty} \frac{2e^{(y-1)z}}{1 + e^{-2z}} = 
		\begin{cases}
			2 &\quad\text{if} ~y =1,\\
			\infty &\quad\text{if} ~y >1.\\
		\end{cases}
		\]
		Hence,
		\[   
		\lim_{z\rightarrow \infty} g(z) =
		\begin{cases}
			\log 2 &\quad\text{if} ~y =1,\\
			\infty &\quad\text{if} ~y >1.\\
		\end{cases}
		\]
		This completes the case $y \ge 1$. Finally, suppose that $y \le -1$. Then, $g'(z) < 0$ for all $z \in \mathbb{R}$, and hence,
		$$\sup_{z\in \mathbb{R}} g(z) = \lim_{z\rightarrow -\infty} g(z)~.$$
		Now, note that
		\[   
		\lim_{z\rightarrow -\infty} e^{g(z)} = \lim_{z\rightarrow -\infty} \frac{2e^{yz}}{e^z + e^{-z}} = \lim_{z\rightarrow -\infty} \frac{2e^{(y+1)z}}{1 + e^{2z}} = 
		\begin{cases}
			2 &\quad\text{if} ~y =-1,\\
			\infty &\quad\text{if} ~y <-1.\\
		\end{cases}
		\]
		Hence,
		\[   
		\lim_{z\rightarrow \infty} g(z) =
		\begin{cases}
			\log 2 &\quad\text{if} ~y =-1,\\
			\infty &\quad\text{if} ~y <-1.\\
		\end{cases}
		\]
		This completes the case $y \le -1$, and the proof of Lemma \ref{sup}. 
	\end{proof}
	
	The following lemma describes the behavior of the function $H_{\beta,p}$.
	
	\begin{lemma}\label{hpr}
		Suppose that $\beta>\beta^*(p)$. Then, the following are true.
		\begin{enumerate}
			\item $H_{\beta,2}' > 0$ on $(0,m_*(\beta,2))$ and $H_{\beta,2}' < 0$ on $(m_*(\beta,2),1)$.
			\item If $p\ge 3$, then $H_{\beta,p}$ can have at most $3$ positive stationary points. Further, there exists $\um(\beta,p) \in (0,m_*(\beta,p))$ such that $H_{\beta,p}' \le 0$ on $(0,\um(\beta,p))$, $H_{\beta,p}' \ge 0$ on $(\um(\beta,p),m_*(\beta,p))$ and $H_{\beta,p}' \le 0$ on $(m_*(\beta,p),1)$.
		\end{enumerate}
	\end{lemma}
	
	\begin{proof}
		To begin with, define $N_{\beta,p}(x) := (1-x^2)H_{\beta_p}''(x)$. Then, 
		$$N_{\beta,p}'(x) = \beta p(p-1)x^{p-3}(p-2-px^2)~.$$ Let us first consider the case $p\ge 3$. Since $N_{\beta,p}'$ has exactly $1$ root in $(0,1)$, it follows by repeated applications of Rolle's theorem, that $H_{\beta,p}'$ can have at most $3$ roots in $(0,1)$. Define:
		$$\um(\beta,p):= \sup\{t\in(0,1]: H_{\beta,p}'\le 0~\textrm{on}~(0,t]\}~.$$
		Since $H_{\beta,p}'(x) = \beta p x^{p-1} - \tanh^{-1}(x)$ and $\lim_{x\rightarrow 0} \tanh^{-1}(x)/x = 1$, we have $\um(\beta,p) > 0$. Clearly, $H_{\beta,p}' \le 0$ on $(0,\um(\beta,p)]$. On the other hand, since $m_*(\beta,p)$ is a global maximizer of $H_{\beta,p}$, and since $\hbp$ can have at most finitely many stationary points, we must have $H_{\beta,p}'(x) >0$ for some $x < m_*(\beta,p)$. This shows that $\um(\beta,p) < m_*(\beta,p)$. Now, by definition of $\um(\beta,p)$, there must exist a sequence $x_n \downarrow \um(\beta,p)$, such that $H_{\beta,p}'(x_n) > 0$ and $x_n > \um(\beta,p)$ for all $n$. Continuity of $H_{\beta,p}'$ now implies that $\um(\beta,p)$ is a stationary point of $\hbp$.   
		
		We will now show that $H_{\beta,p}' \ge 0$ on $(\um(\beta,p),m_*(\beta,p))$. Suppose towards a contradiction, that $H_{\beta,p}'(y) < 0$ for some $y \in (\um(\beta,p),m_*(\beta,p))$. Then, there exist $y_1 \in (\um(\beta,p),y)$ and $y_2 \in (y,m_*(\beta,p))$, such that $H_{\beta,p}'(y_1) >0$ and $H_{\beta,p}'(y_2)>0$. This creates two extra stationary points of $\hbp$, one within $(y_1,y)$ and the other within $(y,y_2)$, giving a total of at least $4$ positive stationary points of $\hbp$, a contradiction! Hence,  $H_{\beta,p}' \ge 0$ on $(\um(\beta,p),m_*(\beta,p))$.
		
		Finally, we show that $H_{\beta,p}' \le 0$ on $(m_*(\beta,p),1)$. Once again, suppose towards a contradiction, that $H_{\beta,p}'(y) > 0$ for some $y \in (m_*(\beta,p),1)$. Since $m_*(\beta,p)$ is a global maximizer of $H_{\beta,p}$ and $\lim_{x\rightarrow 1} H_{\beta,p}'(x) = -\infty$, there exist $y_1 \in (m_*(\beta,p),y)$ and $y_2 \in (y,1)$, such that $H_{\beta,p}'(y_1) <0$ and $H_{\beta,p}'(y_2)<0$. This creates two extra stationary points of $\hbp$, one within $(y_1,y)$ and the other within $(y,y_2)$, giving a total of at least $4$ positive stationary points of $\hbp$, a contradiction! Hence,  $H_{\beta,p}' \le 0$ on $(m_*(\beta,p),1)$. This completes the proof of part (2) of Lemma \ref{hpr}.
		
		Now, suppose that $p=2$. Since $N_{\beta,2}'$ has exactly one root in $(-1,1)$, it follows by repeated applications of Rolle's theorem, that $H_{\beta,2}'$ can have at most $3$ roots in $(-1,1)$. Since $H_{\beta,2}'$ is an odd function, it follows that it can have at most $1$ positive root, which must be $m_*(\beta,2)$. Hence, $H_{\beta,2}'$ must be non-zero and cannot change sign on each of the intervals $(0,m_*(\beta,2))$ and $(m_*(\beta,2),1)$. Since $m_*(\beta,2)$ is a global maximizer of $H_{\beta,2}$, we must have $H_{\beta,2}'(x) >0$ for some $x \in (0,m_*(\beta,2))$, and since $\lim_{x\rightarrow 1} H_{\beta,2}'(x) = -\infty$, we must have $H_{\beta,2}'(x)<0$ for some $x \in (m_*(\beta,2),1)$. Hence, we must have $H_{\beta,2}'(x) >0$ for all $x \in (0,m_*(\beta,2))$ and $H_{\beta,2}'(x)<0$ for all $x \in (m_*(\beta,2),1)$. This proves (1), and completes the proof of Lemma \ref{hpr}.
	\end{proof}

	\begin{lemma}\label{mprop}
		The function $\xi_p(\beta) := m_*(\beta,p)$ is continuous and strictly increasing on $(\beta^*(p),\infty)$. Further, $$\lim_{\beta\rightarrow \infty} \xi_p(\beta) = 1~.$$
	\end{lemma}
	\begin{proof}
		Fix $\beta \in (\beta^*(p),\infty)$ and take a sequence $\beta_n \rightarrow \beta$. Then, $\beta_n \in (\beta^*(p),\infty)$ for all large $n$, and hence, $H_{\beta_n,p}$ will have a unique global maximizer $m_*(\beta_n,p) \in (0,1)$ for all large $n$. Take a subsequence $\{n_k\}_{k\ge 1}$ of the positive integers. This subsequence must have a further subsequence $\{n_{k_\ell}\}_{\ell \ge 1}$ such that $m_*(\beta_{n_{k_\ell}},p) \rightarrow m'$ for some $m' \in [0,1]$. Clearly, $H_{\beta_{n_{k_\ell}},~p}(m_*(\beta_{n_{k_\ell}},p)) \rightarrow H_{\beta,p}(m')$. Since $H_{\beta_{n_{k_\ell}},~p}(m_*(\beta_{n_{k_\ell}},p)) \ge H_{\beta_{n_{k_\ell}},~p}(x)$ for all $x \in [0,1]$ and for all large $\ell$, we must have $H_{\beta,p}(m') \ge H_{\beta,p}(x)$ for all $x \in [0,1]$ (taking $\lim_{\ell \rightarrow \infty}$ on both sides). This means that $m'$ is a non-negative global maximizer of $H_{\beta,p}$. Since $m_*(\beta,p)$ is the only non-negative global maximizer of $H_{\beta,p}$, it follows that $m'=m_*(\beta,p)$. Hence, $m_*(\beta_{n_{k_\ell}},p) \rightarrow m_*(\beta,p)$, showing that $\xi_p(\beta_n) \rightarrow \xi_p(\beta)$, and thereby establishing continuity of $\xi_p$.
		
		Next, take any $t\in (0,1)$, whence $H_{\beta,p}'(t) = \beta p t^{p-1} - \tanh^{-1}(t) >0$ for all $\beta$ large enough. On the other hand, it follows from Lemma \ref{hpr}, that $H_{\beta,p}' \le 0$ on $[m_*(\beta,p),1]$. This shows that $m_*(\beta,p)> t$ for all $\beta$ large enough, showing that
		$\lim_{\beta\rightarrow \infty} \xi_p(\beta) = 1$.
		
		Finally, to show that $\xi_p$ is increasing on $(\beta^*(p),\infty)$, take $\beta_2>\beta_1>\beta^*(p)$. Then, by Lemma \ref{hpr}, $H_{\beta_2,p}' \le 0$ on $[m_*(\beta_2,p),1]$. Since $H_{\beta_1,p}' < H_{\beta_2,p}'$ on $(0,1]$, we must have $H_{\beta_1,p}' < 0$ on $[m_*(\beta_2,p),1]$. However, since $m_*(\beta_1,p)$ is a global maximizer of $H_{\beta_1,p}$, and since $H_{\beta_1,p}$ can have at most finitely many stationary points, there must exist $\varepsilon > 0$, such that $H_{\beta_1,p}' > 0$ on $(m_*(\beta_1,p)-\varepsilon,m_*(\beta_1,p))$. Continuity of $H_{\beta_1,p}'$ now implies that $m_*(\beta_2,p) > m_*(\beta_1,p)$, proving that $\xi_p$ is strictly increasing. This completes the proof of Lemma \ref{mprop}.
	\end{proof}
	
	\section{Technical Results for the Hypergraph Erd\H{o}s-R\'enyi Ising Model}\label{ertech}
	In this section, we prove some technical results related to the hypergraph \ern~Ising model. We start with the proof of Lemma \ref{hamiltcomp}.
	
	\subsection{Proof of Lemma \ref{hamiltcomp}}
	To begin with, for every $\bt \in \sa$, let us define the set:
	$$\Lambda_n(\bt) := \{(i_1,\ldots,i_p)\in [n]^p: x_{i_1}\ldots x_{i_p} =1\}~.$$ Also, let $\ml := \sum_{(i_1,\ldots,i_p) \in \Lambda_n(\bt)} A_{i_1\ldots i_p}$. In these notations, we have:
	$$H_n(\bt) = \alpha_n^{-1} n^{1-p} \left(2\ml - \sum_{(i_1,\ldots,i_p) \in [n]^p} A_{i_1\ldots i_p}\right)~.$$
	For each $\gamma > 0$, define an event:
	$$\Omega_n(\gamma) := \left\{\sup_{\bt \in \sa}\left|\frac{\ml}{\e \ml} - 1\right|  \leq \gamma \right\}~.$$
	Since $\sup_{\bt \in \sa} |\e \ml| \leq \alpha_n n^p$, we have the following on the event $\Omega_n(\gamma_n)$:
	\begin{eqnarray}
		&&\frac{1}{n}\sup_{\bm x \in \sa} |H_n(\bt) - \e H_n(\bt)| \nonumber\\&\leq& 2\alpha_n^{-1}n^{-p} \sup_{\bt \in \sa}|\ml-\e \ml| + \left|\alpha_n^{-1} n^{-p} \sum_{(i_1,\ldots,i_p) \in [n]^p} A_{i_1\ldots i_p} -1 \right|\nonumber\\&\leq& 2\gamma_n ~+~ \Big|\alpha_n^{-1} n^{-p} \sum_{(i_1,\ldots,i_p) \in [n]^p} A_{i_1\ldots i_p} -1 \Big|\label{decomp1}~.
	\end{eqnarray}
	It follows from Theorem 4 in \cite{chernoff_bin}, that 
	\begin{equation}\label{secpart}
		\p\left(\Big|\alpha_n^{-1} n^{-p} \sum_{(i_1,\ldots,i_p) \in [n]^p} A_{i_1\ldots i_p} -1 \Big| > \gamma_n\right) \leq 2e^{-\frac{1}{3}\gamma_n^2\alpha_n n^p} = 2e^{-3n}
	\end{equation}
	In view of \eqref{decomp1}, \eqref{secpart} and the Borel-Cantelli Lemma, it thus suffices to show that
	\begin{equation}\label{toshow}
		\p\left(\Omega_n(\gamma_n)~\textrm{occurs for all but finitely many}~n\right) = 1~,
	\end{equation}
	in order to complete the proof of Lemma \ref{hamiltcomp}. Towards this, note that by a union bound,
	\begin{equation}\label{unioncompl}
		\p(\Omega_n(\gamma)^c) \leq \sum_{\bt \in \sa}\p\left(\ml> (1+\gamma) \e \ml  \right)+ \sum_{\bt \in \sa}\p\left(\ml< (1-\gamma) \e \ml  \right)~.
	\end{equation}
	It follows from Theorem 1 in \cite{hoeffding}, that
	\begin{equation}\label{expmark1}
		\p\left(\ml> (1+\gamma) \e \ml  \right) = \p\left(\frac{\ml}{|\Lambda_n(\bt)| }> (1+\gamma)\alpha_n  \right) \leq e^{-|\Lambda_n(\bt)| D((1+\gamma)\alpha_n\| \alpha_n) }~,
	\end{equation}
	and
	\begin{equation}\label{expmark2}
		\p\left(\ml< (1-\gamma) \e \ml  \right) = \p\left(\frac{\ml}{|\Lambda_n(\bt)| }< (1-\gamma)\alpha_n  \right) \leq e^{-|\Lambda_n(\bt)| D((1-\gamma)\alpha_n\| \alpha_n) }~,
	\end{equation}
	where $D(x\|y) := x\log\frac{x}{y} + (1-x)\log\left(\frac{1-x}{1-y}\right)$. Also, let $$\mn := \left\{-1, -1+\frac{2}{n}, \ldots, 1-\frac{2}{n}, 1\right\}$$ denote the set of all values $\oss := n^{-1} \sum_{i=1}^n x_i$ can take, for some $\bt \in \sa$. Combining \eqref{unioncompl}, \eqref{expmark1} and \eqref{expmark2}, we have by Lemma \ref{card1} and Equation (2.17) in \cite{bovierjsp},
	\begin{eqnarray}\label{uc}
		&&\p(\Omega_n(\gamma)^c)\nonumber\\ &\leq& \sum_{\bt \in \sa} \left\{ e^{-|\Lambda_n(\bt)| D((1+\gamma)\alpha_n\| \alpha_n) }+ e^{-|\Lambda_n(\bt)| D((1-\gamma)\alpha_n\| \alpha_n) } \right\}\nonumber\\&=& \sum_{m \in \mathcal{M}_n} \binom{n}{n(1+m)/2} \left\{e^{-\frac{1}{2} n^p(1+m^p) D((1+\gamma)\alpha_n\| \alpha_n)} + e^{-\frac{1}{2} n^p(1+m^p) D((1-\gamma)\alpha_n\| \alpha_n)} \right\}\nonumber\\&=& e^{-n\left[\frac{n^{p-1}}{2} D((1+\gamma)\alpha_n\| \alpha_n) - \log 2\right] -\frac{\log n}{2} + O(n^{-1})} \sum_{m \in \mn} e^{-\frac{n^p m^p}{2} D((1+\gamma)\alpha_n\| \alpha_n) - n I(m)}\nonumber\\&+& e^{-n\left[\frac{n^{p-1}}{2} D((1-\gamma)\alpha_n\| \alpha_n) - \log 2\right] -\frac{\log n}{2} + O(n^{-1})} \sum_{m \in \mn} e^{-\frac{n^p m^p}{2} D((1-\gamma)\alpha_n\| \alpha_n) - n I(m)}~,
	\end{eqnarray}
	where $I(t) := \frac{1}{2} (1+t)\log(1+t) + \frac{1}{2}(1-t)\log(1-t)$. Since the functions $D$ and $I$ are non-negative, we have
	$$\sum_{m \in \mn} e^{-\frac{n^p m^p}{2} D((1 \pm \gamma)\alpha_n\| \alpha_n) - n I(m)} \leq |\mn| = n+1~.$$
	Hence, we have from \eqref{uc} and Equation (2.30) in \cite{bovierjsp},
	\begin{eqnarray}\label{final}
		\p(\Omega_n(\gamma_n)^c) &\leq& O(\sqrt{n}) \left(e^{-n\left[\frac{n^{p-1}}{2} D((1+\gamma_n)\alpha_n\| \alpha_n) - \log 2\right]} + e^{-n\left[\frac{n^{p-1}}{2} D((1-\gamma_n)\alpha_n\| \alpha_n) - \log 2\right]}\right)\nonumber\\ &\leq& O(\sqrt{n}) \left(\exp\left\{-n\left[\alpha_n n^{p-1}\frac{\gamma_n^2}{6} - \log 2\right]\right\} + \exp\left\{-n\left[\alpha_n n^{p-1}\frac{\gamma_n^2}{4} - \log 2\right]\right\}\right)\nonumber\\&\leq& O(\sqrt{n}) \exp(-0.8 n)~.
	\end{eqnarray}
	Since $\sum_{n =1}^\infty \p(\Omega_n(\gamma_n)^c) < \infty$, \eqref{toshow} follows from \eqref{final} and the Borel-Cantelli lemma, completing the proof of Lemma \ref{hamiltcomp}. \qed
	
	\subsection{Proof of Corollary \ref{hcm1}}
	For each $1\le i\le n$, define $A_{i_2\ldots i_p}^{(i)} := A_{ii_2\ldots i_p}$. Then, for each $1\le i\le n$, one can view $m_i(\bm X)$ as the Hamiltonian (scaled by $n^{-1}$) of the $(p-1)$-spin Erd\H{o}s-R\'enyi Ising model with adjacency tensor $\bm A^{(i)} := ((A_{i_2\ldots i_p}^{(i)}))$. The rest of the proof will follow exactly as the proof of Lemma \ref{hamiltcomp}.

	\subsection{Other Technical Lemmas}
	In this section, we prove some other technical lemmas required for the proof of the main results in Section \ref{sec:ermodel}. We begin with deriving the cardinality of $\Lambda_n(\bt)$ that was required in the proof of Lemma \ref{hamiltcomp}. 
	
	\begin{lemma}\label{card1}
		For every $\bt \in \sa$, we have:
		$$|\Lambda_n(\bt)| = \frac{1}{2}n^p(1+\oss^p)~,$$ where $\oss:= n^{-1}\sum_{i=1}^n x_i$.
	\end{lemma}
	\begin{proof}
		First, note that $(i_1,\ldots,i_p) \in \Lambda_n(\bt)$ if and only if $x_{i_\ell}=-1$ for an even number of $\ell \in [p]:=\{1,\ldots,p\}$. Now, it is easy to see that the number of indices $i \in [n]$ for which $x_i = -1$, is given by $n(1-\oss)/2$. To form an $(i_1,\ldots,i_p) \in \Lambda_n(\bt)$, we must thus choose an even number of these $p$ indices from the total number of $n(1-\oss)/2$ possible indices where we have $-1$, and the rest of these $p$ indices from the remaining $n(1+\oss)/2$ number of possible indices where we have $+1$. We thus have:
		\begin{eqnarray}\label{binexp}
			|\Lambda_n(\bt)| &=& \sum_{k\in [p]\bigcup \{0\}:~k~\textrm{is even}} \binom{p}{k} \left(\frac{n(1-\oss)}{2}\right)^k \left(\frac{n(1+\oss)}{2}\right)^{p-k}\nonumber\\&=& \frac{1}{2}\left(\frac{n(1+\oss)}{2} + \frac{n(1-\oss)}{2}\right)^p + \frac{1}{2}\left(\frac{n(1+\oss)}{2} - \frac{n(1-\oss)}{2}\right)^p~.
		\end{eqnarray}
		Lemma \ref{card1} follows from \eqref{binexp}.
	\end{proof}
	
	The following lemma is crucial in showing consistency of the MPL estimator in the hypergraph Erd\H{o}s-R\'enyi Ising model.
	
	\begin{lemma}\label{c11}
		With probability $1$, we have the following: $$\max_{1\le i_1\le n}\sum_{(i_2,\ldots,i_p) \in [n]^{p-1}} A_{i_1\ldots i_p} = O\left(\alpha_n n^{p-1}\right)~.$$
	\end{lemma}
	\begin{proof}
		Note that $\sum_{(i_2,\ldots,i_p) \in [n]^{p-1}} A_{i_1\ldots i_p} \sim \mathrm{Bin}(n^{p-1},\alpha_n)$. So, by Theorem 4 in \cite{chernoff_bin}, we have:
		$$\p\left(\alpha_n^{-1}n^{1-p}\sum_{(i_2,\ldots,i_p) \in [n]^{p-1}} A_{i_1\ldots i_p} \ge 1+\delta\right) \le e^{-\frac{\delta^2}{2+\delta}\alpha_n n^{p-1}}$$ for every $\delta > 0$. Hence,
		$$\p\left(\max_{1\le i_1\le n} \sum_{(i_2,\ldots,i_p) \in [n]^{p-1}} A_{i_1\ldots i_p} \ge 2(1+\delta) \alpha_n n^{p-1}\right) \le n e^{-\frac{\delta^2}{2+\delta}\alpha_n n^{p-1}} = e^{\log n -\frac{\delta^2}{2+\delta}\alpha_n n^{p-1}}~.$$
		Since $\alpha_n =\Omega(n^{1-p} \log n)$, we can choose $\delta>0$ large enough, so that $\log n -\frac{\delta^2}{2+\delta}\alpha_n n^{p-1} \le -2\log n$ thereby ensuring that
		$$\p\left(\max_{1\le i_1\le n} \sum_{(i_2,\ldots,i_p) \in [n]^{p-1}} A_{i_1\ldots i_p} \ge 2(1+\delta) \alpha_n n^{p-1}\right) \le n^{-2}~.$$ It now follows by an application of the Borel-Cantelli lemma, that:
		$$\p\left(\max_{1\le i_1\le n} \sum_{(i_2,\ldots,i_p) \in [n]^{p-1}} A_{i_1\ldots i_p} \le 2(1+\delta) \alpha_n n^{p-1}~\text{for all large n}\right) = 1~,$$ which completes the proof of Lemma \ref{c11}. 
	\end{proof}
	
\end{document}